\documentclass[11pt,reqno]{amsart}
\usepackage{fullpage}
\usepackage[mathscr]{eucal}
\usepackage{mathrsfs}
\usepackage{amsfonts}
\usepackage{amsmath}
\usepackage{amsthm}
\usepackage{amssymb}
\usepackage{latexsym}
\usepackage[all]{xy}
\usepackage{pstricks,pst-node}
\usepackage[mathscr]{eucal}
\usepackage{palatino}
\usepackage{fullpage}

\theoremstyle{plain}
\newtheorem{thm}[equation]{Theorem}
\newtheorem{cor}[equation]{Corollary}
\newtheorem{prop}[equation]{Proposition}
\newtheorem{lem}[equation]{Lemma}

\theoremstyle{definition}
\newtheorem{defn}[equation]{Definition}

\theoremstyle{remark}

\newtheorem{rem}[equation]{Remark}
\newtheorem{claim}[equation]{Claim}
\newtheorem{notation}[equation]{Notation}

\makeatletter
\renewcommand{\subsection}{\@startsection{subsection}{2}{0pt}{-3ex
plus -1ex minus -0.2ex}{-2mm plus -0pt minus
-2pt}{\normalfont\bfseries}} \makeatother

\numberwithin{equation}{subsection}

\newcommand{\br}[1]{{\mbox{$\left(#1\right)$}}}

\newcommand{\Lmod}[1]{#1\text{-}{\mathsf{mod}}}

\newcommand{\hdot}{{\:\raisebox{2pt}{\text{\circle*{1.5}}}}}
\newcommand{\idot}{{\:\raisebox{2pt}{\text{\circle*{1.5}}}}}

\DeclareMathOperator{\res}{{\mathrm{res}}}

\DeclareMathOperator{\Tor}{\mathrm{Tor}}

\DeclareMathOperator{\sym}{\mathrm{Sym}}

\DeclareMathOperator{\im}{\mathrm{Im}}
\DeclareMathOperator{\supp}{\mathrm{Supp}}
\DeclareMathOperator{\coh}{{\mathrm{Coh}}}

\DeclareMathOperator{\End}{\mathrm{End}}

\DeclareMathOperator{\gr}{\mathrm{gr}}
\DeclareMathOperator{\hilb}{{\mathrm{Hilb}}}

\DeclareMathOperator{\rad}{\mathrm{rad}}

\DeclareMathOperator{\Lie}{\mathrm{Lie}}

\DeclareMathOperator{\ggr}{{\widetilde{\mathrm{gr}}}}
\DeclareMathOperator{\dgr}{{\widetilde{\mathrm{dgr}}}}
\DeclareMathOperator{\Ad}{\mathrm{Ad}}
\DeclareMathOperator{\ad}{\mathrm{ad}}

\DeclareMathOperator{\mo}{{\,\mathrm{mod}\,[\b,\b]}}

\newcommand{\dis}{\displaystyle}

\newcommand{\beq}{\begin{equation}\label}
\newcommand{\eeq}{\end{equation}}

\DeclareMathOperator{\Spec}{\mathrm{Spec}}
\DeclareMathOperator{\pr}{pr}

\newcommand{\iso}{{\;\stackrel{_\sim}{\to}\;}}
\newcommand{\cd}{\!\cdot\!}
\DeclareMathOperator{\Hom}{\mathrm{Hom}}
\DeclareMathOperator{\hhom}{{{\scr H}\!{\mathit{om}}}}

\def\ccirc{{{}_{\,{}^{^\circ}}}}
\newcommand{\bo}{\mbox{$\bigotimes$}}
\renewcommand{\o}{\otimes }

\newcommand{\cc}{{\scr C}}
\newcommand{\kap}{{\kappa}}
\newcommand{\kkap}{{\wt\kappa}}
\newcommand{\si}{\sigma }
\newcommand{\cchi}{{\boldsymbol{\chi}}}

\renewcommand{\t}{{\mathfrak t}}

\newcommand{\wt}{\widetilde }
\newcommand{\UU}{{\mathfrak U}}

\newcommand{\eu}{{\operatorname{\mathsf{s}}}}
\newcommand{\eps}{{\epsilon}}
\newcommand{\ttr}{{{\mathfrak T}^r}}
\newcommand{\tts}{{\overset{_{_\circ}}{\mathfrak T}}}
\newcommand{\ts}{{\overset{_{_\circ}}{\mathfrak t}}}
\newcommand{\Id}{{\operatorname{Id}}}
\newcommand{\zz}{{\mathfrak C}}
\newcommand{\xx}{{\mathfrak X}}
\newcommand{\xc}{{\mathfrak X}^\circ}

\newcommand{\fl}{{\mathfrak l}}

\newcommand{\ha}{{\mathfrak W}}
\newcommand{\yy}{{\mathfrak x}}
\newcommand{\tx}{{\wt{\mathfrak X}}}
\newcommand{\zr}{{\mathfrak C}^{rs}}

\newcommand{\xr}{{\mathfrak X}^{rs}}
\newcommand{\erem}{\hfill$\lozenge$\end{rem}\vskip 3pt }

\newcommand{\scr}[1]{\mathscr{#1}}
\newcommand{\lo}{\stackrel{L}{\mbox{$\bigotimes$}}}
\newcommand{\loo}{{\stackrel{_L}\otimes}}

\newcommand{\mhm}{Hodge module }
\newcommand{\HH}{{\mathbb H}}
\newcommand{\fa}{{\mathfrak a}}

\def\ccirc{{{}_{^{\,^\circ}}}}

\newcommand{\DD}{{\boldsymbol{\mathcal D}}}
\renewcommand{\ss}{{\mathcal{S}}}
\newcommand{\gl}{{\mathfrak g\mathfrak l}}

\newcommand{\filt}{F^{\operatorname{ord}}}
\newcommand{\ord}{{\operatorname{ord}}}
\newcommand{\hodge}{{{\operatorname{Hodge}}}}

\newcommand{\dd}{{\mathscr{D}}}
\newcommand{\oo}{{\mathcal{O}}}

\renewcommand{\aa}{{\mathscr{A}}}

\newcommand{\mh}{{\mathsf{HM}}}
\renewcommand{\tt}{{\mathfrak T}}
\newcommand{\h}{{\mathfrak t}}
\newcommand{\hh}{{{\mathfrak T}}}

\newcommand{\dcoh}{D^b_{\operatorname{coh}}}

\renewcommand{\H}{{\scr H}}

\newcommand{\NN}{{\mathcal N}}
\newcommand{\tg}{{\wt\g}}

\renewcommand{\gg}{{\mathfrak G}}

\newcommand{\hc}{{Harish-Chandra }}

\newcommand{\TT}{{T}}
\newcommand{\I}{{\mathcal I}}
\newcommand{\J}{{\mathcal J}}

\newcommand{\C}{\mathbb{C}}
\newcommand{\g}{\mathfrak{g}}
\renewcommand{\b}{\mathfrak{b}}

\newcommand{\norm}{{\operatorname{norm}}}
\renewcommand{\red}{{\operatorname{red}}}
\newcommand{\La}{\Lambda }

\newcommand{\pp}{{\scr P}}
\newcommand{\rr}{{\scr R}}

\newcommand{\bb}{{\scr B}}

\newcommand{\D}{{\mathsf D}}
\newcommand{\x}{{\mathfrak x}}
\newcommand{\mm}{{\mathcal M}}

\newcommand{\T}{{\mathcal T}}
\newcommand{\inv}{^{-1}}

\newcommand{\Z}{{\mathbb Z}}
\newcommand{\en}{{\enspace}}

\newcommand{\vi}{${\en\sf {(i)}}\;$}
\newcommand{\vii}{${\;\sf {(ii)}}\;$}
\newcommand{\viii}{${\sf {(iii)}}\;$}
\newcommand{\iv}{${\sf {(iv)}}\;$}

\newcommand{\sset}{\subset}
\newcommand{\sminus}{\smallsetminus}
\newcommand{\into}{\,\,\hookrightarrow\,\,}
\newcommand{\too}{\,\,\longrightarrow\,\,}
\newcommand{\mto}{\mapsto}
\newcommand{\onto}{\,\,\twoheadrightarrow\,\,}
\newcommand{\N}{{\mathcal{N}}}
\newcommand{\tgg}{{\wt{\mathfrak G}}}
\newcommand{\mmu}{{\boldsymbol{\mu}}}
\newcommand{\nnu}{{\boldsymbol{\nu}}}
\newcommand{\ppi}{{\boldsymbol{\pi}}}

\newcommand{\Ga}{\Gamma }

\newcommand{\om}{\omega }

\renewcommand{\sl}{{\mathfrak s\mathfrak l}}
\begin{document}
\title{{\large{\textsf{Isospectral  commuting variety and\\
the Harish-Chandra {\large$\boldsymbol{\mathcal D}$}-module}}}}

\author{Victor Ginzburg}\address{
Department of Mathematics, University of Chicago,  Chicago, IL 
60637, USA.}
\email{ginzburg@math.uchicago.edu}

\begin{abstract} Let  $\g$ be  a complex reductive Lie algebra 
with Cartan algebra $\h$. 
Hotta and Kashiwara defined a
holonomic $\dd$-module $\mm$, on $\g\times\h$, called
Harish-Chandra module. We  give an explicit
description of 
$\gr\mm$, the associated
graded module with respect to
 a canonical  {\em Hodge filtration} on $\mm$.
The description involves
the {\em isospectral
commuting variety}, a subvariety $\xx\sset\g\times\g\times\h\times\h$
which is a finite extension of the variety of pairs of commuting
elements of $\g$.
Our main result establishes an isomorphism
of $\gr\mm$ with  the structure
sheaf of the normalization of $\xx$.
It follows, thanks to the theory of polarized
Hodge modules, that
the  normalization of the isospectral
commuting variety is    Cohen-Macaulay
and Gorenstein, confirming a conjecture of M. Haiman.

In the special case where $\g=\gl_n$,
there is an open subset  of the
  isospectral
commuting variety that is closely related to 
the Hilbert scheme of $n$ points in $\C^2$.
The sheaf $\gr\mm$
gives rise to a locally free sheaf on the Hilbert scheme.
 We show that the corresponding vector bundle is
isomorphic to the 
{\em Procesi bundle} that plays an important role 
in the work of M. Haiman.
\end{abstract}
\maketitle


\section{The isospectral commuting variety}\label{main_sec}
\subsection{Reminder on commuting schemes}
Let $G$ be a connected complex reductive group with
Lie algebra $\g$. 
We fix $T\sset G$, a maximal torus, and
let $\t=\Lie T$ be the corresponding Cartan subalgebra 
of $\g$.
The group $G$ acts
on $\g$ via the adjoint action $G\ni g: x\mto \Ad g(x).$

We put
$\gg:=\g\times\g$ and let $G$ act diagonally on $\gg$.
The {\em commuting scheme} $\zz$, of the Lie
algebra $\g$, is defined as the {\em scheme-theoretic}
zero fiber of the commutator map $\kap:\ \gg\to\g,$ $(x,y)\mto[x,y]$.
Thus, $\zz$ is  a closed 
$G$-stable  subscheme of 
$\gg$ and, set-theoretically, one has
$\zz=\{(x,y)\in\gg\mid [x,y]=0\}$.

Given $x\in \g$, we let $G_x=\{g\in G\mid \Ad g(x)=x\}$
and write $\g_x=\Lie G_x$ for the centralizer of $x$ in $\g$.
The isotropy group of a pair  $(x,y)\in\gg$ under the $G$-diagonal
action
equals $G_{x,y}:= G_x\cap G_y$.
Let $\g_{x,y}:=\Lie G_{x,y}=\g_x\cap\g_y$ denote
the corresponding Lie algebra.
We call an element $x\in\g$, resp.  $(x,y)\in\zz$,
{\em regular} if we have $\dim\g_x=\dim\t$,
resp.  $\dim\g_{x,y}=\dim\t.$
Write $\g^r$, resp. $\zz^r$, for the set of regular elements
of $\g$, resp. of $\zz$.
 Further, let  $\zz^{rs}\sset\zz^r$ be the set
of  pairs $(x,y)\in\zz$ such that both $x$ and $y$ are
regular  semisimple elements.

Basic properties of the commuting scheme  may be summarized
as follows.
\begin{prop}\label{zzbasic} \vi The set $\zz^{rs}$ is 
 Zariski open and dense in $\zz$.

\vii The smooth
locus of the scheme $\zz$ equals $\zz^r$.
\end{prop}

Here, part (i) is due  to Richardson \cite{Ri1}.
To prove (ii), one views the commutator map 
as a moment map $\kap: T^*\g\to\g^*$.
The result then follows easily  from a general
property of   moment maps 
saying that $\ker(d\kap)$,
the kernel of the
differential of a moment map $\kap$, equals the annihilator
of the image of $(d\kap)^\top$,
 the dual map,
cf.  also \cite[Lemma 2.3]{Po}.
\medskip

Proposition \ref{zzbasic} implies that 
$\zz$ is a generically reduced and irreducible 
scheme.
It is a long standing open problem whether or not
this scheme is reduced.
\medskip

Given a scheme
$X$,  write $\oo_X$  for the structure
sheaf and put $\C[X]=\Ga(X,\oo_X)$.
Let $X_\red$ denote the  reduced scheme corresponding
to $X$.
Given an algebraic
group $K$ acting on $X$,  let $\C[X]^K\sset\C[X]$ denote
the subalgebra of $K$-invariants,
and put $X/\!/K:=\Spec \C[X]^K$.

Let $N(T)$ be the normalizer of $T$ in $G$ and 
$W=N(T)/T$ be the Weyl group.
We put $\tt:=\t\times\t$ and let $W$ act diagonally on $\tt$.
Clearly, $\tt\sset\zz$. So,   restriction of polynomial functions
gives  algebra maps 
\beq{Jo} 
\res:\ \C[\gg]^G\, \onto\,\C[\zz]^G\,\to\,\C[\tt]^W.
\eeq

The group $W$ acts freely on $\tt^r:=\zz^r\cap\tt$,
a Zariski open dense subset of $\tt$.
Further, the assignment
$(g,x,y)\mto \big(\Ad g(x),\,\Ad g(y)\big)$ 
 induces an isomorphism
$G\times_{N(T)}\ttr\iso\zz^{rs}.$
It follows, since the set  $\zz^{rs}$ is dense in
$\zz$ and 
the second  map  in \eqref{Jo} is surjective
by a theorem  of Joseph \cite{Jo},
that this map induces an isomorphism
$\tt/W\iso[\zz/\!/G]_\red.$ 
That isomorphism is analogous to the
isomorphism $\t/W\iso \g/\!/G$ induced
by the Chevalley isomorphism
$\C[\g]^G\iso\C[\t]^W$.
Thus, we have
a diagram 
\beq{X}
\xymatrix{
\tt\ \ar@{->>}[r]&\ \tt/W\
 \ar@{=}[r]&\ [\zz/\!/G]_\red\
\ar@{^{(}->}[r]&\ \zz/\!/G\ &\ \zz.\ar@{->>}[l]
}
\eeq

In the special case where
$\g=\gl_n$, the scheme $\zz/\!/G$ is known
to be reduced,  by \cite[Theorem 1.3]{GG}.
It  is expected to be reduced 
for any reductive
Lie algebra $\g$.

\subsection{}
There is a natural $G\times W$-action
on $\g\times\t$,
 resp. on $\gg\times\tt$,
induced by the   $G$-action  on the first factor and 
 $W$-action on the second factor.
There is also 
a $\C^\times$-action on $\g\times\t$
by dilations,
and an induced $\C^\times\times\C^\times$-action
on $\gg\times\tt=(\g\times\t)\times(\g\times\t)$.

The fiber product $\yy:=\g\times_{\g/\!/G}\h,$ resp. 
 $\zz\times_{\zz/\!/G}\tt$, is
a closed $G\times W\times\C^\times$-stable
subscheme in $\g\times\t,$ resp. 
$G\times W\times\C^\times\times\C^\times$-stable
subscheme in $\gg\times\tt$.
The first projection $\yy\to\g,$ resp.
$\zz\times_{\zz/\!/G}\tt\to\zz$,
is a $G$-equivariant {\em finite} morphism.
The group
$W$ acts along the fibers of this morphism.

\begin{lem}\label{xxbasic} \vi The set
$\yy^{rs}:=\g^{rs}\times_{\g/\!/G}\h^r$,
resp.
$\xr:=\zz^{rs}\times_{\zz/\!/G}\tt^r$,
 is  a smooth, irreducible, Zariski open and dense
subset of $\yy$, resp. of $\zz\times_{\zz/\!/G}\tt$.

\vii The first projection  $\yy^{rs}\to\g^{rs},$ resp.
$\xr\onto \zr$, is
a Galois covering
with Galois group ~$W$.
\end{lem}

For the proof that $\xr$ is  irreducible and Zariski  dense
in  $\zz\times_{\zz/\!/G}\tt$ see Corollary \ref{xdense}
of \S\ref{END} below. All other statements of the
above lemma are clear.

The scheme $\yy$ is known to be a reduced normal
complete intersection in $\g\times\t$,
cf. eg \cite{BB}.
On the contrary, the scheme $\zz\times_{\zz/\!/G}\tt$
is not  reduced
already in the case $\g=\sl_2$.

\begin{defn}\label{xxdef} The
{\em isospectral commuting variety}
is defined as
$\dis\xx:=[\zz\times_{\zz/\!/G}\tt]_\red,$
a  {\em reduced} fiber product. Let
 $p_{_\zz}:\ \xx\to\zz$, resp. $p_{_\tt}:\ \xx\to\tt$, 
denote the first, resp. second, projection.
\end{defn}
\begin{rem}
The isospectral commuting  variety has been considered by M. Haiman,
in \cite[\S8]{Ha2}, \cite[\S 7.2]{Ha3}.
\end{rem}

Lemma \ref{xxbasic}
shows that $\xx$ is an irredicible variety and that
we may (and will) identify
the set $\xr$ with a Zariski open subset of $\xx$.

Let $q: T^*X\to X$ denote the cotangent bundle on
a smooth variety $X$. The group $\C^\times$ acts along the fibers of $q$
by dilations. An 
irreducible reduced subvariety $\Lambda \sset T^*X$
is said to be a {\em Lagrangian  cone} if 
it is $\C^\times$-stable and the tangent space to $\La$ at any smooth
point of $\La$
is a Lagrangian vector
subspace with respect to the
canonical symplectic 2-form on ~$T^*X.$

We use an invariant
bilinear form on $\g$ to identify
$\g^*$ with $\g$, resp. $\t^*$ with
$\t$.
This gives an identification $T^*(\g\times\t)=\gg\times\tt$.
The $\C^\times$-action along the fibers of $q$
corresponds to the action of the subgroup $\{1\}\times\C^\times
\sset\C^\times\times\C^\times$.

Let $\yy^r:=\{(x,t)\in \yy\mid x\in \g^r\}$.
This is  a Zariski open and dense subset
of $\yy$ which is
contained in the smooth locus of $\yy$
(since the differential of the adjoint quotient
map $\g\to\g/\!/G$ has maximal rank at any point of
$\g^r$).
Let $N_{\x^r}\sset T^*(\g\times\t)$
be the total space of the conormal
bundle on $\x^r$. Thus,
$\overline{N_{\x^r}}$, the closure of $N_{\x^r}$,
is a Lagrangian cone in $T^*(\g\times\t)$.

\begin{lem}\label{conormal} In $T^*(\g\times\t)=\gg\times\tt$,
we have $\overline{N_{\x^r}}=\xx$; 
in particular, $\xx$ is a Lagrangian cone.
\end{lem}
\begin{proof} We know that $\xx$ is 
a reduced and irreducible, $\C^\times$-stable
subvariety. By Lemma \ref{xxbasic}, one has
$\dim\xx=\dim\xr=\dim\zz^{rs}
=\dim(G\times_{N(T)}\tt^r)=
\dim\g+\dim\t$.
Hence,  $\dim\xx=\frac{1}{2}\dim T^*(\g\times\t)$.
We leave to the reader to check
that  tangent spaces to $\xr$ are  isotropic subspaces
 with respect to
the symplectic form on $T^*(\g\times\t)$.
It follows that $\xx$ is 
 a  Lagrangian cone.
Further,  set theoretically, one has $q(\xx)=\g\times_{\g/\!/G}\t
=\x$.

The above properties force $\xx$ to be equal to the closure
of the total space of the conormal bundle on the smooth
locus of $\x$,
cf. eg. \cite{CG}, Lemma 1.3.27.
\end{proof}

\subsection{An analogue of the Grothendieck-Springer
resolution}\label{tx}
We recall a few  definitions.

Let $\bb$ be the flag variety, the variety
of all Borel subalgebras $\b\sset\g$.
For any pair $\b,\b'$, of  Borel subalgebras, there
is a canonical isomorphism $\b/[\b,\b]\cong\b'/[\b',\b']$,
cf. eg. \cite[Lemma 3.1.26]{CG}.
This quotient 
is referred to as the `abstract
Cartan algebra'.

Given a  Borel subgroup $T\sset B\sset G$,
 we may identify $\bb\cong G/B$.
Write $\b=\Lie B$, resp. $\t=\Lie T$.
The composite map $\t\into\b\onto \b/[\b,\b]$
induces an isomorphism $\t\iso\b/[\b,\b]$.
So, we will abuse the notation and write
$\t$ for the abstract Cartan algebra as well.

Motivated by Grothendieck and Springer,
we introduce the following varieties
$$
\tg:=\{(\b,x)\in\bb\times\g\mid x\in\b\},\en\text{resp.}\en
\tgg:=\{(\b,x,y)\in\bb\times\g\times\g\mid
x,y\in\b\}.$$

The first projection
  makes $\tg$, resp. $\tgg$,
a sub vector bundle of the trivial vector bundle
$\bb\times\g\to\bb$, resp. $\bb\times\gg\to\bb$. One has 
 a $G$-equivariant vector bundle isomorphism 
$\tg\cong G\times_{B}\b,$ resp. $\tgg\cong G\times_{B}(\b\times\b)$.
Thus, $\tg$ and $\tgg$ are smooth connected varieties.
Further, the  assignment $(\b,x)\mto x\mo\in\b/[\b,\b],$
resp. $(\b,x,y)\mto  (x\mo,\ y\mo),$
gives
a well defined smooth morphism $\nu:\ \tg\to\h$,
resp. $\nnu:\ \tgg\to\hh$.
There is also a projective morphism
 $\mu:\ \tg\to\g,\ 
(\b,x)\mto x$,
resp. $\mmu:\ \tgg\to \gg,\ (\b,x,y)\mto(x,y)$.
The map $\mu$ is known as  the
{\em Grothendieck-Springer resolution}.

\begin{notation}
Let $\T_X$ denote the tangent sheaf, resp. $K_X:=\wedge^{\dim X} \T^*_X$ 
denote the canonical sheaf (or line bundle), on a smooth variety $X$.
\end{notation}

\begin{prop}\label{tb}
\vi The image of the map $\mu\times\nu:\tg\to\g\times\t$ is
contained in $\yy=\g\times_{\gg/\!/G}\t$.
The resulting morphism
$\pi:\tg\to\yy$ is a  resolution 
of singularities such that $\pi:\ \pi\inv(\yy^r)\iso\yy^r$ is an
isomorphism.

\vii The image of the map $\mmu\times\nnu:
\tgg\to\gg\times\tt$ is
contained in $\gg\times_{\gg/\!/G}\tt$.
The resulting map gives a resolution
of singularities 
 $\tgg\onto[\gg\times_{\gg/\!/G}\tt]_\red$.

\viii The canonical bundle on $\tg$ has a natural trivialization
by a nowhere vanishing section $\omega\in K_\tg$ such that one has
$\mu^*(dx)=(\prod_{\alpha>0}\ \langle \alpha,\nu\rangle)\cdot
\omega$, where the product in the RHS is taken over 
all positive roots  $\alpha\in\t^*$ and where $dx$ is a constant
volume form on $\g$.

\iv We have $K_{\tgg}=\wedge^{\dim\bb}\ {\mathfrak q}^* \T_\bb,$
where ${\mathfrak q}:\ \tgg\to\bb$ stands for the first projection.
\end{prop}

Part (i) of this proposition is well-known,
cf. \cite{BB}, \cite{CG}.
Part (ii) is an immediate consequence
of Lemma \ref{lss}, of \S\ref{ss1} below.
This lemma implies, in particular, that the image of the
map $\mmu$ equals the set of pairs $(x,y)\in\gg$ such that
$x$ and $y$ generate a solvable Lie subalgebra of $\g$.
The statement in (ii) is a variation of
 results concerning the
null-fiber of the adjoint
quotient map $\gg\to\gg/\!/G$, see
\cite{Ri2}, \cite{KW}.

The 
descriptions of canonical bundles in 
Proposition \ref{tb}(iii)-(iv) are straightforward.
Finally, the equation of part (iii) that involves
the section $\omega$
appears eg. in  \cite{HK1}, formula (4.1.4).
This equation is, in essence, nothing but Weyl's 
classical integration formula.\qed
\medskip

Let $\N$ be  the variety of nilpotent elements
of $\g$ and put
$\wt\N:=
\{(\b,x)\in\bb\times\g\mid x\in[\b,\b]\}$.
One may restrict the commutator map $\kap$
to each Borel subalgebra $\b\sset\g$ to obtain
 a morphism
$\kkap:\ \tgg\to \wt\N,\
(\b,x,y)\mto (\b,\,[x,y])$.
One has a commutative diagram,
where the vertical map in the middle is known as the Springer resolution,
\beq{kk}
\xymatrix{
\bb=\mu\inv(0)\ \ar[d]^<>(0.5){\mu} \ar@{^{(}->}[rr]^<>(0.5){\imath:
\b\mto(\b,0)}&&\ 
\wt\N\
\ar@{->>}[d]^<>(0.4){\mu|_{_{\wt\N}}}&&\ \tgg\ \ar[ll]_<>(0.5){\kkap}
\ar[d]^<>(0.5){\mmu}\\
\{0\}\  \ar@{^{(}->}[rr]&&
\ \NN\ \ar@{^{(}->}[r]&\ \g\ &\ \gg\ \ar[l]_<>(0.5){\kap}
}
\eeq

We introduce the following closed subset in $\tgg$:
\beq{txform}\tx :=
\{(\b,x,y)\in\bb\times\g\times\g\mid
x,y\in \b,\en [x,y]=0\}\ =\ \kkap\inv\big(\imath(\bb)\big).
\eeq
Diagram \eqref{kk} shows that the morphism $\mmu$ maps $\tx$ to $\zz$,
resp.
 $\mmu\times\nnu$ maps $\tx$ to
$\zz\times_{\zz/\!/G}\hh.$

\begin{rem} Let
$\tx^{rs}:=\tx\,\cap\, \nnu\inv(\tt^r).$
This is an open subset of $\tx$.
One can show that $\tx^{rs}$ is a smooth local complete intersection in 
$\tgg$ and the map $\mmu\times\nnu$ induces 
 an isomorphism
$\tx^{rs}\iso\xx^{rs}$.
It is likely that the scheme $\tx$ is not
irreducible so that $\tx^{rs}$
 is not  dense in $\tx$.
\end{rem}

\subsection{A DG algebra}\label{dima_sec}
Let  $\T:=\T_\bb$ and write $q: \wt\N\to \bb$,
resp. ${\mathfrak q}: \tgg\to\bb$, for the first   projection.
The cotangent space at a point $\b\in\bb$ may
be identified with the vector space
$(\g/\b)^*\cong [\b,\b]$. This yields 
 a natural isomorphism $\wt\N\cong T^*\bb$,
of $G$-equivariant vector bundles on $\bb$,
cf. \cite[ch. 3]{CG}.
For each $n\geq 0,$ we put $\aa^n=
\wedge^n\,{\mathfrak q}^*\T$,
 a coherent  sheaf on $\tgg$.
The sheaf ${\mathfrak q}_*\aa^n$ is a
 quasi-coherent sheaf on $\bb$ that may be identified
with the  sheaf of sections of
a $G$-equivariant vector bundle on
$\bb$ with fiber
$\sym (\b^*\oplus\b^*)\o\wedge^n [\b,\b]^*.$

The sheaf $q^*\T^*$, on $T^*\bb$,
has  a canonical section  $\eu$.
Using the isomorphisms
$T^*\bb=\wt\N$ and $\kkap^*(q^*\T^*)={\mathfrak q}^*\T^*$,
we may view $\kkap^*\eu$
as a section  of the sheaf ${\mathfrak q}^*\T^*$.
Contraction with $\kkap^*\eu$
gives a morphism $i_{\kkap^*\eu}:\
\wedge^\hdot\,{\mathfrak q}^*\T\to
\wedge^{\hdot-1}\,{\mathfrak q}^*\T$.
This makes
$\aa^\hdot:=
(\wedge^\hdot\,{\mathfrak q}^*\T^*,\,i_{\kkap^*\eu})$
a DG algebra, with differential
$i_{\kkap^*\eu}$ of degree $(-1)$.

Now, view the commutator map  as an element
$\kap\in\b^*\o\b^*\o[\b,\b]=
\Hom(\b\o\b,\ [\b,\b])$.
The differential
on the DG algebra ${\mathfrak q}_\idot\aa^\hdot$,
induced by $i_{\kkap^*\eu}$, acts 
by  fiberwise contraction with the  element
$\kap$
as follows (for any $k,m,n\geq 0$):
\beq{bundle}
i_{\kap}:\
\sym^k\b^*\o\sym^m\b^*\o(\wedge^n [\b,\b]^*)
\to\sym^{k+1}\b^*\o\sym^{m+1}\b^*\o(\wedge^{n-1} [\b,\b]^*).
\eeq

The DG algebra $\aa^\hdot$ has an important 
 self-duality property. Observe first that, in degree
  $d:=\dim \bb$,  the top
nonzero degree of the algebra $\aa^\hdot$,
we have $\aa^d=K_\tgg$, by Proposition \ref{tb}(iv).
Moreover, multiplication in $\aa^\hdot$ induces
a canonical isomorphism of complexes
\beq{selfdual}
\aa^{d-\hdot}\iso\hhom_{_{\oo_{\tgg}}}(\aa^\hdot, \aa^d)=
 \hhom_{_{\oo_{\tgg}}}(\aa^\hdot, K_{\tgg}).
\eeq

\begin{notation}
Let $\dcoh( X)$ denote the bounded derived category
of the category $\coh X$, of coherent sheaves on
a scheme $X$. Given ${\mathcal F}\in \dcoh( X)$,
let $\H^j({\mathcal F})\in\coh X$
denote the $j$th cohomology sheaf, resp.
 $\HH^j(X, {\mathcal F})$
denote the $j$th hyper-cohomology group, of ${\mathcal F}$.
 \end{notation}

Associated with diagram \eqref{kk},
there are natural derived functors
\beq{ddd}
 \xymatrix{\dcoh( \bb)\
\ar[r]^<>(0.5){\imath_*}&\ 
\dcoh( \wt\N)\ar[r]^<>(0.5){L\kkap^*}\
&\ \dcoh( \tgg)\ \ar[rr]^<>(0.5){R(\mmu\times\nnu)_*}
&&\ \dcoh(\gg\times\tt).}
\eeq

Contraction with the canonical section $\eu$ yields
a Koszul complex $\ldots\to  \wedge^3q^*\T\to\wedge^2q^*\T\to
q^*\T\to \imath_*\oo_\bb\to0$.
This is the standard  locally
free resolution of the sheaf $\imath_*\oo_\bb$, on $\wt\N$.
Therefore, the DG algebra $\aa^\hdot\cong\kkap^*(\wedge^\hdot\T,\
i_\eu)$ provides a
DG algebra model for the object $L\kkap^*(\imath_*\oo_\bb)
\in \dcoh(\tgg)$.
The  corresponding DG scheme $\Spec\aa^\hdot$
may be thought of as a `derived
analogue' of the scheme $\wt\xx$, in the sense
of derived algebraic geometry. 
The coordinate ring of the DG scheme $\Spec\aa^\hdot$
is  $R\Ga(\tgg,\aa^\hdot)=R\Ga(\bb,{\mathfrak q}_*\aa^\hdot)$,
a DG algebra 
well defined up to quasi-isomorphism, cf. \cite{Hi}.
Hence, $\HH^\hdot(\tgg,\aa^\hdot)$,
the hyper-cohomology of that DG algebra,
acquires the canonical structure of a
graded commutative
algebra. 
The following result says  that the
 DG scheme $\Spec\aa^\hdot$ is, in a sense, a `resolution'
of the isospectral commuting variety.

\begin{thm}\label{dima} We have
$\HH^k(\tgg, \aa^\hdot)=0$ for all $k\neq 0$,
and there is a canonical $G\times \C^\times\times\C^\times$-equivariant
algebra isomorphism
$\HH^0(\tgg, \aa^\hdot)\cong\C[\xx_\norm]$,
where $\xx_\norm$ is the normalization of 
 the isospectral commuting variety.
\end{thm}

The statement of Theorem \ref{dima} (in an equivalent form
presented in Theorem  \ref{dima2} below) was
 suggested to me by Dmitry Arinkin.

\subsection{Bigraded character formula}
Given a reductive group $K$, a $K$-scheme $X$ and a $K$-equivariant
coherent sheaf $\mathcal F$, on $X$, we write $\cchi^K(\mathcal F)$
for its equivariant
Euler characteristic,
a class in the Grothendieck group of rational
$K$-modules.
We will identify the Grothendieck group of  rational
$\C^\times\times\C^\times$-modules
with $\C[q_1^{\pm 1}, q_2^{\pm 1}]$,
a Laurent polynomial ring.

Recall that we have a natural
$G\times W\times \C^\times\times\C^\times$-action
on the variety
$\xx$.
Since a  reductive group action can always 
be lifted canonically to
the normalization,
cf. \cite{Kr}, \S4.4, we obtain
a $G\times W\times \C^\times\times\C^\times$-action
on  $\xx_\norm$ as well.
This makes the coordinate ring
$\C[\xx_\norm]$ a bigraded $G\times W$-module.

Let $R_+\sset \t^*$ denote the
set of weights of the adjoint $\t$-action on the vector space
$[\b,\b]^*\cong\g/\b$,
and let $\ell(-)$ denote the length function
on the Weyl group $W$.

\begin{cor}\label{character_cor} The 
bigraded $T$-character of the coordinate ring of the variety $\xx_\norm$
is given by the formula
$$
\cchi^{T\times \C^\times\times\C^\times}
(\oo_{\xx_\norm})=
\sum_{w\in W}  (-1)^{\ell(w)}\!\!
 \prod_{\alpha,\beta_1,\beta_2,\gamma\in R_+}
\frac{(1-q_{_1}q_{_2}e^{w(\alpha)})}{(1-q_{_1}e^{w(\beta_1)})(1-q_{_2}
e^{w(\beta_2)})(1-e^{-w(\gamma)})}.
$$
\end{cor}
\begin{proof} For any locally finite representation
$E$, of a Borel subgroup $B\sset G$, 
let $\underline{E}$ denote the
corresponding induced $G$-equivariant vector bundle
on $\bb=G/B$.

We let the group $\C^\times$ act
on $\b$ by dilations. This gives
a $\C^\times\times\C^\times$-action along
the fibers, $\b\times\b$, of the projection ${\mathfrak q}:\
\tgg\to\bb$. That makes $\tgg$
a $G\times \C^\times\times\C^\times$-variety.
For each $m\geq 0$, the sheaf 
$\aa^m$ comes equipped with a natural
$G\times \C^\times\times\C^\times$-equivariant structure.
From the definition of the sheaf
${\mathfrak q}_\idot\aa^\hdot$, we get an equation, cf. \eqref{bundle}:
\beq{above}
\cchi^{G\times \C^\times\times\C^\times}({\mathfrak q}_\idot\aa^\hdot)=
\sum_{k,m,n\geq 0}\ (-q_1q_2)^n
q_1^{k}q_2^{m}\cdot \cchi^G\big(
\sym^k\underline{\b}^*\o\sym^m\underline{\b}^*\o(\wedge^n 
\underline{[\b,\b]}^*)\big).
\eeq

Thanks to the
$G\times\C^\times\times\C^\times$-equivariant
 isomorphism of  Theorem \ref{dima},
the LHS of  equation \eqref{above} is equal to
 the bigraded character of the $G$-module $\C[\xx_\norm]$.
To obtain the formula of the corollary,
one computes
the  RHS of  equation 
\eqref{above} 
using the Atiyah-Bott fixed point formula
for $G$-equivariant vector bundles on the flag variety,
cf. eg. \cite[\S 6.1.16]{CG}.
\end{proof}

It would be very interesting to describe the
structure of $\C[\xx_\norm]$ as a $W$-module
in terms of the DG algebra $\aa^\hdot$.

\subsection{Main results}\label{dsec}
Given an irreducible,  not necessarily reduced, scheme $X$
let $\psi: X_\norm\to X_\red$ denote the normalization of
the corresponding reduced scheme $X_\red$.

One may view the complex $(\aa^\hdot, i_{\kkap^*\eu})$
as an object of $\dcoh(\tgg)$. Thus, taking the (derived) 
direct image via the morphism
$\mmu\times\nnu$, we get
\beq{global}
R\Ga(\tgg,\,\aa^\hdot)=
R\Ga\big(\gg\times\tt,\ R(\mmu\times\nnu)_*\aa^\hdot\big)=
R\Ga\big(\gg\times\tt,\ R(\mmu\times\nnu)_*L\kkap^*(\imath_*\oo_\bb)\big).
\eeq

Using that the scheme $\gg\times\tt$ is affine, 
we see from the above isomorphisms that Theorem \ref{dima}
is  equivalent to the following sheaf theoretic result

\begin{thm}\label{dima2} The 
sheaves $\H^k\big(R(\mmu\times\nnu)_*L\kkap^*(\imath_*\oo_\bb)\big)$
vanish for all $k\neq 0$ and
there is an isomorphism
$\H^0\big(R(\mmu\times\nnu)_\idot\aa^\hdot\big)\cong\psi_\idot\oo_{\xx_\norm}$,
of sheaves of $\oo_{\gg\times\tt}$-algebras.
\end{thm}

This theorem will be deduced
from  Theorem \ref{cohiso}  of \S\ref{filtmmG}
below.

The significance of Theorem \ref{dima2}
is due to the self-duality isomorphism \eqref{selfdual}.
It follows,  since $\mmu\times\nnu$
 is a proper morphism,
 that the object
$R(\mmu\times\nnu)_*\aa^\hdot\in\dcoh(\gg\times\tt)$ is
isomorphic to its Grothendieck-Serre dual, up to a shift.
Therefore, the vanishing
statement in Theorem  \ref{dima2} 
implies, in particular, that the sheaf 
$\H^0\big(R(\mmu\times\nnu)_\idot\aa^\hdot\big)$ is Cohen-Macaulay.

Thus, Theorem  \ref{dima2} yields
the following result that confirms a conjecture of M. Haiman, \cite[Conjecture
7.2.3]{Ha3}.
\begin{thm}\label{main_thm}  The sheaf $\oo_{\xx_\norm}$ is Cohen-Macaulay
and $\xx_\norm$ is a Gorenstein variety with trivial
canonical bundle.\qed
\end{thm}

The composite map $p_{_\zz}\ccirc\psi:\ \xx_\norm\to\xx\to\zz$
factors through  a finite 
$G$-equivariant morphism $p:\ \xx_\norm\to\zz_\norm.$
We put  $\rr:=p_\idot\oo_{\xx_\norm}$,
a $G$-equivariant coherent sheaf of $\oo_{\zz_\norm}$-algebras.
The group $W$ acts along the fibers of $p$ and
this gives a $W$-action on $\rr$
by $\oo_{\zz_\norm}$-algebra automorphisms.
Write $\rr^W$ for the subsheaf of
$W$-invariant sections.

From Theorem \ref{main_thm} we deduce

\begin{prop}\label{vectbun} \vi The  canonical
morphism $\oo_{\zz_\norm}\to\rr^W$
is an isomorphism, equivalently, the
map $p$ induces an isomorphism
$\xx_\norm/W\iso \zz_\norm$.

\vii  The sheaves
$\rr$ and $\rr^W$ are
Cohen-Macaulay; thus,
$\zz_\norm$ is a Cohen-Macaulay variety.

\viii
The restriction of the sheaf $\rr$
to  $\zz^r,$ the smooth locus of $\zz$,
 is a locally free sheaf.
Each fiber of the corresponding algebraic vector bundle
affords
the regular representation of 
the group ~$W$.
\end{prop}

\begin{proof}  The scheme $\xx_\norm/W$ is reduced
and integrally
closed,
as a quotient of an integrally
closed  reduced scheme
by a finite group action, \cite{Kr}, \S3.3. 
Further, by  Lemma
\ref{xxbasic}(ii), the map   $p$ is generically a Galois
covering with $W$ being the Galois group.
It follows that the induced map $\xx/W\to\zz_\norm$
is finite and birational. Hence it is
  an isomorphism and  (i) is proved.

Part (ii) follows from  Theorem \ref{main_thm}
since the property of a coherent sheaf  be
Cohen-Macaulay is stable under taking  direct 
images by finite
morphisms and also under taking direct
summands.
Finally,
any coherent Cohen-Macaulay sheaf  on a smooth variety is 
locally free, hence $\rr|_{\zz^r}$ is a locally free
sheaf. 
 Lemma
\ref{xxbasic}(ii) implies that
the fiber of the corresponding vector bundle over any point 
of $\zr$ affords
the regular representation of the group
$W$. 
By continuity, the same holds
for the fiber at any  point of $\zz^r$, since
 $\zr$ is dense in ~$\zz$.
\end{proof}

\section{The \hc $\dd$-module}
\subsection{Algebraic definition of the \hc module}\label{r_sec}
Let $\dd_X$ denote the sheaf of 
 algebraic differential
operators on a smooth algebraic variety $X$,
and write $\dd(X):=\Ga(X,\dd_X)$ for the algebra of
global sections. Given  an action on $X$ of an algebraic
group $K$, we let  $\dd(X)^K\sset \dd(X)$
denote
 the subalgebra of $K$-invariant 
differential operators.


We have
 $\dd(\g)$ and $\dd(\t)$, the algebras of polynomial differential operators
on the vector spaces $\g$ and $\t$,
respectively. The corresponding subalgebras
 of  differential operators  with
constant coefficients may be identified with
$\sym\g$ and $\sym\t,$ the symmetric algebras
of $\g$ and $\t$, respectively.
Further, given $a\in \g$, one may view
the  map $\ad a:\ \g\to\g,\ x\mto [a,x]$ as a (linear)  vector
field on $\g$, that is, 
 as a first order differential operator on $\g$.
The assignment $a\mto\ad a$ gives
a linear map $\ad: \g\to\dd(\g)$, with image $\ad\g$. 

\hc  defined a `radial
part' map
\beq{r}
\rad:\ \dd(\g)^G\to\dd(\t)^W,
\eeq
an algebra homomorphism
such that its restriction to $G$-invariant polynomials,
resp.  to  $G$-invariant
constant coefficient  differential operators,
reduces to
the Chevalley isomorphism
$\C[\g]^G\iso \C[\t]^W$,
resp. $(\sym\g)^G\iso (\sym\t)^W$.


\begin{rem} According to the results of Wallach \cite{Wa} and
Levasseur-Stafford \cite{LS1}-\cite{LS2},
the radial part map \eqref{r} induces an algebra
isomorphism
$$
\rad:\ \dd(\g)^G/[\dd(\g)\ad\g]^G=
\big[\dd(\g)/\dd(\g)\ad\g\big]^G\iso \dd(\t)^W.
\eqno\lozenge
$$
\end{rem}

We will use
a special notation $\DD:=\dd_{\g\times\t}$
for the sheaf of  differential operators
on the vector space $\g\times\t$.
We have $\Ga(\g\times\t, \DD)=\dd(\g\times\t)=\dd(\g)\o\dd(\t)$.

\begin{defn}\label{mm_def} 
The \hc module is a left $\DD$-module defined as follows
\beq{M}
\mm\
:=\ \DD\big/\br{\DD\cd(\ad\g\o 1)+\DD\cd\{u\o1-1\o \rad(u),\ u\in \dd(\g)^G\}}.
\eeq
\end{defn}

\begin{rem} The above definition was motivated by,
but is not identical to, the definition
of  Hotta and Kashiwara, see \cite{HK1},
formula (4.5.1). The equivalence of the two definitions
follows from Lemma \ref{I=I}, of \S4 below.
\end{rem}

\subsection{The order filtration on the \hc module}\label{hc_filt}
Let $X$ be a smooth variety.
The  sheaf $\dd_X$ comes equipped with
an ascending filtration $\filt_\idot\dd_X$, by the order of
 differential operator. 
For the associated graded sheaf,
one has a canonical isomorphism $\gr^\ord\dd_X\cong
q_\idot\oo_{T^*X},$
where $q: T^*X\to X$ is the cotangent bundle projection.

Let  $M$ be a $\dd_X$-module.
An
ascending filtration
 $F_\idot M$ such that one has $\filt_i\dd_X\cd F_j M\sset
F_{i+j}M,\ \forall i,j,$
is said to be {\em good} if
$\gr^F M$, the associated graded module,
is a coherent $q_\idot \oo_{T^*X}$-module.
In that case, there is a canonically defined
 coherent sheaf $\ggr^F M$,
on $T^*X$, such that one has an
isomorphism $\gr^F M=q_\idot\ggr^F M$,
of   $q_\idot \oo_{T^*X}$-modules.
We write $[\supp(\ggr^F M)]$ for the {\em support cycle}
of the sheaf $\ggr^F M$. This  is  an algebraic
cycle in $T^*X$, a linear
combination of the irreducible components
of  the support of $\ggr^F M$ counted with multiplicites.
The cycle  $[\supp(\ggr^F M)]$ is known to be independent
of the choice of a good filtration on $M$, cf.
\cite{Bo}.

Recall that  $M$ is called {\em holonomic}
if one has $\dim \supp(\ggr^F M)=\dim X$.
In such a case,  each irreducible
component of $\supp(\ggr^F M)$, viewed as a reduced
variety, is a Lagrangian
cone in $T^*X$.
Holonomic $\dd_X$-modules form an abelian category
and the assignment $M\mto [\supp(\ggr^F M)]$
is additive on short exact sequences of holonomic modules,
cf. eg. \cite{Bo},  \cite{HTT}.
From this, one obtains

\begin{lem}\label{CC}  If $M$ is a holonomic
$\dd_X$-module such that
the cycle $[\supp(\ggr^F M)]$ equals the fundamental cycle of
an irreducible variety taken with multiplicity 1,
then $M$ is a simple  $\dd_X$-module.
\end{lem}
\medskip

According to formula \eqref{M} 
 the \hc module has the form $\mm=\DD/\I$,
where  $\I$ is a left ideal of $\DD$.
The order filtration on  $\DD$ restricts to a filtration
on $\I$ and it also induces a quotient filtration 
$\filt_\idot\mm$, on $\DD/\I$.
Using the identifications $T^*(\g\times\t)=
\gg\times\tt$ and
 $\ggr^\ord\DD=\oo_{\gg\times\tt},$ we have
$\ggr^\ord\mm=\oo_{\gg\times\tt}/\ggr^\ord\I$,
where $\ggr^\ord\I$, the associated graded
ideal, is a subsheaf of ideals of $\oo_{\gg\times\tt}$.

Let  $\J\sset \oo_{\gg\times\tt}$ be the ideal sheaf
of  $\zz\times_{\zz/\!/G}\tt$,
a (non reduced) subscheme in
$\gg\times\tt$.
The following  result provides a relation between the ideals
$\ggr^\ord\I$ and $\J$.

\begin{lem}\label{supp}  One has inclusions
$\J\sset\ggr^\ord \I\sset\sqrt{\J}$, of ideals, and
an equality $[\supp(\ggr^\ord\mm)]=[\xx]$,
of algebraic cycles in $\gg\times\tt$.
In particular,
 $\mm$ is a nonzero simple holonomic $\DD$-module.
\end{lem}
\begin{proof} 
To simplify notation, it will be convenient below
to  work with spaces of global
sections rather than with sheaves.
Thus, let $\D:=\dd(\g\times\t)=\dd(\g)\o\dd(\t)$ be
 the  algebra
of polynomial differential operators
on  $\g\times\t$ so, one has $\gr^\ord\D=\C[\gg\times\tt]$. Put
$I=\Ga(\g\times \t, \I)$,
resp. $J=\Ga(\gg\times\tt, \J)$.
The  vector space $\g\times\t$ being affine, 
it is sufficient to prove an analogue 
of the lemma for the ideals $\gr^\ord I$ and $J$.

Let $\kap^*: \g=\g^*\to\C[\gg]$
be the pull-back morphism
induced by the commutator map $\kap$.
By definition, one has
$\C[\zz]=\C[\gg]/\C[\gg]\kap^*(\g)$.
Hence, we obtain
$$
\C[\zz\times_{\zz/\!/G}\ \tt]=
\C[\zz]\bo_{\C[\zz]^G}\C[\tt]=
\big(\C[\gg]/\C[\gg]\kap^*(\g)\big)\ \bo_{\C[\gg]^G}\ \C[\tt].
$$
We conclude that
$J$
is an ideal of the algebra $\C[\gg\times\tt]=\C[\gg]\o\C[\tt]$
 generated by the following set
$$
\overline{S}:=\
\kap^*(\g)\o1\ \cup\ \{f\o1-1\o\res(f)\mid f\in\C[\gg]^G\}.
$$

To compare the ideals $\g^\ord I$ and $J$ we introduce
a subset of the algebra $\D$ as follows:
$$
S:=\ \ad\g\o1\ \cup\
\{u\o1-1\o \rad(u),\ u\in \dd(\g)^G\}.
$$

\begin{claim}\label{claim}
We have  $\overline{S}=\{\gr^\ord(s),\
s\in S\};$ in other words, the set
 $\overline{S}$ is the set
formed by the principal symbols of the  elements
of  $S$.
\end{claim}

To prove the claim, we observe
that the map $\gr^\ord(\ad):\ \g\to \gr^\ord\D=\C[\gg]$
may be identified with the map $\kap^*$
considered above.
Further, the order  filtration  $\filt_\idot\dd(\t)$,
resp. $\filt_\idot\dd(\g)$,
induces a  filtration  on $\dd(\t)^W$,
resp.  on $\dd(\g)^G$.
It follows from definitions that
the radial part map \eqref{r} respects the filtrations and
$\gr^\ord(\rad),$ the associated graded map, is nothing but
the algebra map $\res$ in 
~\eqref{Jo}. This completes the proof of the claim.

Thus, we see that 
principal symbols of   elements
of the set $S$  generate 
 the ideal $J$.
On the other hand, formula
\eqref{M} shows that 
the set $S$ is a set of generators
of the left ideal $I$.
This yields an inclusion
 $J\sset \gr^\ord I$. 
Furthermore, it
follows by a standard argument that,
for  any point
$x\in \gg\times\tt$ where the ideal generated
the principal symbols of  elements 
of the set $S$ is reduced,
one has an equality
$\oo_x\bo_{\C[\gg\times\tt]} J
=\oo_x\bo_{\C[\gg\times\tt]} \gr^\ord I$,
of ideals in the local ring $\oo_x$
of the  point $x$.
In particular, this equality holds
for any point  $x\in \xx^{rs}$,
by Lemma \ref{xxbasic}(i).
Hence,
any element $f\in \gr^\ord I$, viewed as a function
on $\gg\times\tt$,
vanishes on the set $\xx^{rs}.$
The set $\xx^{rs}$
being Zariski dense in $\zz\times_{\zz/\!/G}\ \tt$,
by Lemma \ref{xxbasic}(i),
we deduce that the function $f$ vanishes
on the zero set of  the ideal $J$.
Hence, $f\in \sqrt{J},$
by Hilbert's Nullstellensatz.
Thus, we have proved
that $\gr^\ord I\sset\sqrt{J}.$

Now, the inclusion $J\sset \gr^\ord I$ implies that,
set theoretically, one has
$\supp(\ggr^\ord\mm)\sset\xx$.
Further,   the scheme
$\zz\times_{\zz/\!/G}\ \tt$ being generically
reduced, by Lemma \ref{xxbasic}(i),
 the inclusion $\gr^\ord I\sset\sqrt{J}$
implies that
the fundamental cycle  $[\xx]$ occurs in $[\supp(\ggr^\ord\mm)]$
with multiplicity one.
Finally, the dimension of 
any irreducible component
of the support of the sheaf  $\ggr^\ord\mm$ 
is $\geq\dim\xx,$ since $\xx$ is a Lagrangian subvariety.
We conclude that $\xx$ is the only irreducible
component of $\supp(\ggr^\ord\mm)$.
Thus, we have $[\supp(\ggr^\ord\mm)]=[\xx]$,
and Lemma \ref{CC} completes the proof.
\end{proof}

\subsection{Hodge filtration}\label{hc_hodge}
In his work \cite{Sa}, M. Saito defines,
for any smooth algebraic variety $X$ over $\mathbb Q$,
a semisimple abelian category $\mh(X)$, of polarized
 Hodge modules. The data of a \mhm includes,
in particular,  a holonomic $\dd_X$-module
 $M$, with regular singularities,
and a good filtration $F$, on $M$,
called {\em Hodge filtration}. Abusing 
notation, we write $(M,F)\in \mh(X)$,
and let $\ggr(M,F)$ denote the corresponding
coherent sheaf on $T^*X.$

There is a duality functor
on the abelian category of holonomic $\dd_X$-modules,
called Verdier duality, cf. \cite{HTT}.
One has a similar duality functor ${\mathbb D}:\ \mh(X)\to
\mh(X)$. Saito shows that the assignment $(M,F)\mto \ggr(M,F)$
intertwines the functor  ${\mathbb D}$
with the standard  Grothendieck duality
on $\dcoh( T^*X)$.
Specifically, he proves that, for $(M,F)\in \mh(X)$, the 
Grothendieck dual of the 
coherent sheaf $\ggr(M,F)$ is isomorphic,
up to shifts, to $\ggr({\mathbb D}(M,F))$.
This implies, in particular, that the
Grothendieck dual of  $\ggr(M,F)$, viewed
an object of
$\dcoh( T^*X)$, has a single nonvanishing
cohomology sheaf,
cf.  \S\ref{filt_sec}. In other words, 
one has, see \cite{Sa}, Remark 5.1.4(7).

\begin{lem}\label{cohen}  For any $(M,F)\in \mh(X)$, the
$\oo_{T^*X}$-module $\ggr(M,F)$ is Cohen-Macaulay.
\end{lem}

 Hotta and Kashiwara \cite{HK1} gave
a geometric construction of the \hc module,
to be explained in \S\ref{hk}.
That construction gives
the \hc module $\mm$  the natural structure
of a   polarized Hodge module on $\g\times\t$
which is,
in addition,  isomorphic to its Verdier dual, up to 
a shift.
Thus, there is a canonical
Hodge filtration $F^\hodge\mm$, on $\mm$.
Write $\ggr^\hodge\mm$ for the corresponding
associated graded sheaf.
Applying Lemma
\ref{cohen} and
using the self-duality of $\mm$, we obtain, cf. also Lemma \ref{supp}(ii),

\begin{cor}\label{serre} The
sheaf $\ggr^\hodge\mm$ is a 
Cohen-Macaulay coherent $\oo_{\gg\times\tt}$-module.
This module is
set-theoretically supported on $\xx$ and it
is isomorphic
to its Grothendieck dual, up to a shift.
\end{cor}

In the previous subsection,
we have considered another filtration
$\filt_\idot\mm$,  the order filtration 
on the \hc module. We do not know if the Hodge 
and the order filtrations on $\mm$ are equal or not.
The result below, to be proved in \S\ref{pf},  provides a partial answer.

\begin{lem}\label{ordhodge}
\vi For each $k\geq 0$, one has inclusions
$\jmath:\ F^\ord_k\mm\sset F^\hodge_k\mm$,
and the restriction of the induced map
 $\gr\jmath:\ \ggr^\ord\!\mm\to\ggr^\hodge\mm$
to the open set $\xx^{rs}$ is an isomorphism.

\vii The $\oo_{\gg\times\tt}$-module $\ggr^\hodge\mm$
has a natural
structure of commutative $\oo_{\gg\times\tt}$-algebra
such that  the following composite map  is
a morphism of 
$\oo_{\gg\times\tt}$-algebras
\beq{comp}
\xymatrix{
\oo_{\zz\times_{\zz/\!/G}\ \tt}=\oo_{\gg\times\tt}/\J\
\ar@{->>}[rr]^<>(0.5){\text{Lemma \ref{supp}(i)}}&&\
\oo_{\gg\times\tt}/\I=\ggr^\ord\!\mm\
\ar@{->>}[r]^<>(0.5){\gr\jmath}&\
\ggr^\hodge\mm.
}
\eeq
\end{lem}

\begin{rem} The Hodge filtration on $\mm$ that we
use, in Lemma \ref{ordhodge} in particular,
   differs by a degree shift from the
one used by Saito.  The precise normalization
of our Hodge filtration will be specified
in \S\ref{pf}.
Degree shifts clearly do not
affect the validity of Lemma \ref{cohen},
hence, Corollary \ref{serre} holds with
our normalization of the Hodge filtration
as well.
\end{rem}

\subsection{The Hodge filtration on the \hc module}\label{main_res}
Recall the normalization map  $\psi:\xx_\norm\to\xx$.

One of the main results of the paper
is the following description of the sheaf $\ggr^\hodge\mm$.

\begin{thm}\label{mthm} There is a
natural $\oo_{\gg\times\tt}$-algebra isomorphism
$\psi_\idot\oo_{\xx_\norm}\iso \ggr^\hodge\mm.$
\end{thm}

\begin{rem}
Note that the above theorem, combined with
Corollary \ref{serre}, implies Theorem \ref{main_thm}
from \S\ref{dsec}. This yields a proof of  Theorem \ref{main_thm}
that does not refer to Theorem \ref{dima2}. However,
both Theorem \ref{dima2} and Theorem \ref{mthm}
will be deduced from the same result, Theorem \ref{cohiso},
so these two theorems are closely related.
\hfill$\lozenge$
\end{rem}

\begin{proof}[Sketch of proof of  Theorem \ref{mthm}]
Let  $X:=\Spec_{\gg\times\tt}(\ggr^\hodge\mm)$ be
 the relative spectrum of  $\ggr^\hodge\mm$,
the latter being
viewed, thanks to Lemma \ref{ordhodge}(ii), as a 
 sheaf of commutative $\oo_{\gg\times\tt}$-algebras.
Thus, $X$ is a scheme equipped with a finite
morphism $X\to\gg\times\tt$ that factors through
a morphism $f: X\to \zz\times_{\zz/\!/G}\ \tt$,
by Lemma \ref{ordhodge}(ii).
We know that   ${\zz\times_{\zz/\!/G}\tt}$
 is an irreducible and generically reduced  scheme,
by Lemma \ref{xxbasic}(i).
Furthermore, the scheme $X$  is Cohen-Macaulay,
by Corollary \ref{serre},
and   $f$
is  an isomorphism over the generic point
of $\zz\times_{\zz/\!/G}\ \tt$, by
Lemma \ref{ordhodge}(i).
It follows that $X$ is   reduced.
Therefore, the map $f$ factors through a  finite morphism
$f: X\to\xx$.

Recall next that $\xx$ is irreducible
(Lemma \ref{xxbasic}) and, we have
$\ f_*[X]=[\supp(\ggr^\hodge\mm)]$
$=[\supp(\ggr^\ord\mm)]=[\xx]$,
by  Lemma \ref{supp}.
The equality $f_*[X]=[\xx]$, of algebraic cycles,
 forces $X$ be irreducible and the map
$f: X\to\xx$  be a finite birational isomorphism.

To proceed further
we  need the following technical 
\begin{defn}\label{U}
Let $\zz^{rr}$ be the set of pairs $(x,y)\in\zz$
such that either $x$ or $y$
is a regular element of $\g$.
We put $\xx^{rr}:=p_{_\zz}\inv(\zz^{rr})$, a 
subset of $\xx$.
\end{defn}

The following two lemmas  explain the role of the set $\zz^{rr}$
 in the proof of Theorem  \ref{mthm}.

\begin{lem}\label{irr} \vi $\zz^{rr}$  is a Zariski open subset and
 one has $\dim(\zz\sminus \zz^{rr})\leq \dim\zz-2$;

\vii The set $\xx^{rr}$ is   contained in the smooth locus
of the variety $\xx$ and 
it is  Zariski
dense in $\xx$.
\end{lem}

\begin{rem} \vi  It is clear from Proposition \ref{zzbasic}
that the  set $\zz^{rr}$ 
is contained in the smooth locus
of $\zz$ and that it is  dense in $\zz$.

\vii It is
essential that in  Lemma \ref{irr}(ii)  one takes the preimage
of $\zz^{rr}$ in $\xx$; the preimage of $\zz^{rr}$ in
$\zz\times_{\zz/\!/}\tt$, a {\em nonreduced} fiber product,
is not smooth already for $\g=\sl_2$.

\viii  Definition \ref{U} and
Lemma \ref{irr} were motivated, in
part, by \cite[Lemma 3.6.2]{Ha1}.
\end{rem}

\begin{lem}\label{keylem} The morphism
$f$ induces an isomorphism
$\ f\inv (\xx^{rr})\iso \xx^{rr}$.
\end{lem}

Part (i) of Lemma \ref{irr} will be proved in
\S\ref{ss3} and part (ii)  of Lemma \ref{irr}
is Lemma \ref{codim1}(ii), of \S\ref{ulem} below.
For the proof of Lemma \ref{keylem} see \S\ref{keypf}.

We can now complete the proof of Theorem \ref{mthm}.
The map
 $p_{_\zz}$, resp. $f$,  being finite, from Lemma \ref{irr}(i)
we deduce that the complement of the set $\xx^{rr}$,
resp. of  the set 
$f\inv(\xx^{rr})$, has codimension one in $\xx$, resp. in $X$.
Hence,  Lemma \ref{irr}(ii) implies that the variety $\xx$,
resp. $X$, 
is smooth in codimension 1.
We conclude, since  $X$ is  irreducible and Cohen-Macaulay,
that $X$ is  normal  and, moreover,  
 $f=\psi$ is the normalization map.
\end{proof}

The above proof shows that the isomorphism  of Theorem 
\ref{mthm}
fits into the following chain of  $\oo_{\gg\times\tt}$-algebra maps,
 cf.  Lemma \ref{supp}(ii),
$$
\xymatrix{
\oo_{\zz\times_{\zz/\!/G}\ \tt}\ 
\ar@{->>}[r]&\ \ggr^\ord\!\mm\ \ar@{->>}[r]&
\ \oo_\xx\ \ar@{^{(}->}[r]&
\ \psi_\idot\oo_{\xx_\norm}\ \ar@{=}[rr]^<>(0.5){\text{Theorem \ref{mthm}}}&&
\ \ggr^\hodge\mm.
}
$$
\section{Relation to work of M. Haiman}
Thoughout this section,
we consider a special case where $G=GL_n.$

\subsection{A vector bundle on the Hilbert scheme}\label{psec}
Let
$V=\C^n$ be the fundamental representation of
the group $GL_n$, so $\g=\End_\C V$.
For any $(x,y)\in \zz$, one may view
the vector space $\g_{x,y}$  as an associative
subalgebra of $\End_\C V$.
The group $G_{x,y}$  may be identified
with the group of invertible elements of
that subalgebra, hence, this group is connected.

Let $\C(x,y)$ denote the associative
subalgebra of $\End_\C V$ generated by the elements $x$ and $y$.
Recall that a vector $v\in V$ is said to be 
a {\em cyclic vector} for a pair $(x,y)\in\zz$ 
if one has $\C(x,y)v=V$.
Let $\zz^\circ$ be the set of pairs $(x,y)\in\zz$
which have a cyclic vector.
Part (i) of the following result is due to
Neubauer and Saltman   \cite{NS} and part (ii)
is well-known,  cf.  \cite{NS}.

\begin{lem}\label{cyc} \vi A pair $(x,y)\in\zz$ is regular
if and only if one has
 $\C(x,y)=\g_{x,y}$.

\vii The set  $\zz^\circ$
is a Zariski open subset of $\zz^r$.
For  $(x,y)\in\zz^\circ$,
all cyclic vectors for $(x,y)$ form
a single  $G_{x,y}$-orbit, which is
an open $G_{x,y}$-orbit in $V$. \qed
\end{lem}

Next, introduce a variety of triples
$$\ss=\{(x,y,v)\in \g\times\g\times V\en\big|\en
[x,y]=0,\  \C(x,y)v=V\}.$$

It is easy to see that $\ss$ is a smooth quasi-affine variety.
We let  $G$ act on $\ss$ by
$\,g: (x,y,v)\mto (gxg\inv,gyg\inv,gv)$. This
action  is free, by
Lemma \ref{cyc},
and it is known that
there exists a universal  geometric quotient morphism
$\rho:\ \ss\to\ss/G$. Furthermore, the variety $\ss/G$ 
may be identified with $\hilb:=\hilb^n(\C^2)$, 
the Hilbert scheme of $n$ points in $\C^2$,
 see
\cite{Na}. Writing $\delta$ for 
the projection $(x,y,v)\mto (x,y),$
one obtains
a commutative diagram, cf. also \cite[\S8]{Ha2},
\beq{ss}
\xymatrix{
\zz^\circ\ \ar@{->>}[drr]_<>(0.5){\eqref{X}}
&&\ 
\ss\
\ar@{->>}[ll]_<>(0.5){\delta}\ar@{->>}[rr]^<>(0.5){\rho}&&\ \hilb
\ar@{->>}[dll]^<>(0.5){^{\;\ \text{Hilbert-Chow}}}\\
&& \ \tt/W\ &&
}
\eeq

Let ${\mathcal G}$ be the scheme  theoretic preimage of the
diagonal in $\ss\times\ss$ under the
morphism $G\times\ss\to\ss\times\ss,\ (g,s)\mto (gs,s).$
Set theoretically, one has
${\mathcal G}=\{(x,y,\gamma)\in \g\times\g\times G\mid
{(x,y)\in\zz,}\ \gamma\in G_{x,y}\}$.
The group $G$ acts naturally on ${\mathcal G}$
and the projection ${\mathcal G}\to\zz,
\ (x,y,\gamma)\to (x,y)$ is a $G$-equivariant map.
The scheme ${\mathcal G}$ has the natural
structure of a group scheme over $\zz$.
Lemma \ref{cyc}(i) implies
that ${\mathcal G}|_{\zz^r}$,
the restriction of ${\mathcal G}$ to
the regular locus of the commuting scheme,
is a smooth commutative
group scheme. We put
${\mathcal G}^\circ:={\mathcal G}|_{\zz^\circ}$.

We define a  ${\mathcal G}^\circ$-action along the fibers of the
map $\delta$, see \eqref{ss}, by 
$(x,y,\gamma):\ (x,y,v)\mto (x,y,\gamma v)$.
Lemma \ref{cyc}(ii) implies that
this ${\mathcal G}^\circ$-action on $\ss$ is free.
Moreover,   the map $\delta$ makes $\ss$ a $G$-equivariant
${\mathcal G}^\circ$-torsor over ${\zz^\circ}$. In particular,
$\delta$ is a smooth morphism.

We use the notation of diagram \eqref{ss} and
put $\pp=(\rho_\idot\delta^*\rr)^G$.
By equivariant descent, one has a canonical
isomorphism $\rho^*\pp=\delta^*\rr$.
Proposition
\ref{vectbun} yields the following result

\begin{cor}\label{P} The sheaf $\pp$ 
is a locally free coherent
sheaf of commutative $\oo_{\hilb}$-algebras 
equipped with a natural $W$-action, by algebra automorphisms.
The fibers of 
the corresponding algebraic vector bundle 
afford the regular representation of 
the group ~$W$.\qed
\end{cor}

\subsection{The isospectral Hilbert scheme}
\label{hilb}
Let $\ha:=\Spec_{\hilb}\pp$ be
 the relative spectrum of $\pp$,
a sheaf  of algebras on the Hilbert scheme.
The scheme $\ha$ comes equipped
with a
flat and finite morphism  $\eta: \ha\to\hilb$
and with a $W$-action along the fibers
of $\eta$.
We conclude that
 $\ha$ is a reduced Cohen-Macaulay and Gorenstein variety.

One can interpret the construction of
the scheme $\ha$ in more geometric terms as follows.
Let $\xc:=p\inv({\zz^\circ})$ and consider
the following commutative diagram
\beq{square}
\xymatrix{
\xc\ \ar[dr]_<>(0.5){p}&\ \xc\times_{\zz^\circ}\ss\
\ar@{->>}[l]_<>(0.5){\wt\delta}
\ar@{->>}[dr]^<>(0.5){\wt p}\ar@{.>}[rr]^<>(0.5){h}_<>(0.5){\cong}
&&\ \ss\times_{\hilb}\ha\
\ar[dl]_<>(0.5){\wt\eta}\ar@{->>}[r]^<>(0.5){\wt\rho}
&\ \ha
\ar@{->>}[dl]^<>(0.5){\eta}\\
&\ {\zz^\circ}\ &\ \ss\ \ar@{->>}[l]^<>(0.5){\delta}\ar@{->>}[r]_<>(0.5){\rho}&
\ \hilb&
}
\eeq

In this diagram, the morphisms $\delta$ and $\rho$ are
smooth, resp.  the morphisms $p$ and $\eta$ are finite and flat.
The morphism $\wt\delta$, resp.
$\wt\rho, \wt p,$ and $\wt\eta$, 
is obtained from $\delta$, resp. from $\rho, p,$ and $\eta$,
by base change.
Hence, flat base change yields
\beq{basechange}\delta^*\rr=\delta^*p_\idot\oo_{\xc}=
\wt p_\idot\oo_{_{\xc\times_{\zz^\circ}\ss}},\quad
\text{resp.}\quad
\rho^*\pp=\rho^*\eta_\idot\oo_\ha=
\wt\eta_\idot\oo_{_{\ss\times_{\hilb}\ha}}.
\eeq

We equip the scheme $\ss\times_{\hilb}\ha$ with a $G$-action
induced from the one on $\ss$, resp. equip 
 $\xc\times_{\zz^\circ}\ss$ 
with the
$G$-{\em diagonal} action.
Thus, $\wt p$, resp. $\wt\eta,$ is a
$G$-equivariant finite and flat morphism.
Since $\delta^*\rr=\rho^*\pp$,
from \eqref{basechange}
we deduce  a canonical isomorphism
$\wt p_\idot\oo_{_{\xc\times_{\zz^\circ}\ss}}
\cong\wt\eta_\idot\oo_{_{\ss\times_{\hilb}\ha}}$,
of $G$-equivariant  sheaves of $\oo_\ss$-algebras.
This means that there is a canonical $G$-equivariant  isomorphism
$h:\ \xc\times_{\zz^\circ}\ss\iso
\ss\times_{\hilb}\ha$, the dotted arrow in diagram
\eqref{square}.

The map
$\wt\rho$ makes $\ss\times_{\hilb}\ha$ a
$G$-torsor over $\ha$. Therefore, the composite
$\wt\rho\ccirc h$ makes $\xc\times_{\zz^\circ}\ss$
a $G$-torsor over $\ha$, hence, we have $\ha=
(\xc\times_{\zz^\circ}\ss)/G$, a
 geometric quotient by $G$.
The scheme $\xc$ being Gorenstein,
we deduce that the scheme $\ha$ is Gorenstein as well.

The  $W$-action on $\xc$ induces one on
$\xc\times_{\zz^\circ}\ss$.
This $W$-action commutes with the $G$-diagonal
action, hence, descends
to $(\xc\times_{\zz^\circ}\ss)/G$. 
The resulting  $W$-action may be identified with
the $W$-action on
$\ha$  that was defined earlier.
The composite map
$p_{_\tt}\ccirc\wt\delta:\ \xc\times_{{\zz^\circ}}\ss\to\xc\to\tt$
descends to a $W$-equivariant morphism $\ha\to\tt$.
We obtain the following  diagram
\beq{ha}
\xymatrix{\xc\ &
\ \xc\times_{{\zz^\circ}}\ss\ \ar@{->>}[l]_<>(0.5){\wt\delta}
\ar@{->>}[rr]^<>(0.5){\wt\rho\ccirc h}&&\ 
\  (\xc\times_{{\zz^\circ}}\ss)/G=\ha  \ar[rr]^<>(0.5){\si\times\eta}&&
\tt\times_{\tt/W}\hilb.
}
\eeq

Following Haiman, one defines the
{\em isospectral Hilbert scheme} as
$[\hilb\times_{\tt/W}\tt]_\red,$
a reduced fiber product.
We have

\begin{prop}\label{haim} 
The map $\si\times\eta$, on the right of \eqref{ha},
factors through an isomorphism 
$$\ha\iso [\hilb\times_{\tt/W}\tt]_\norm.$$
In particular, the normalization of the 
isospectral Hilbert scheme is Cohen-Macaulay
and Gorenstein.
\end{prop}
\begin{proof} It is clear 
that $\wt\rho\ccirc h(\wt\delta\inv(\xr))$  is a Zariski
open subset of $\ha$.
Lemma \ref{xxbasic}(ii) implies
readily that the map $\si\times\eta$ restricts to
an isomorphism $\wt\rho\ccirc h(\wt\delta\inv(\xr))$
$\iso\hilb\times_{\tt/W}\ttr$;
in particular, $\si\times\eta$ is a birational isomorphism.
The image of the map $\si\times\eta$
 is contained
in $[\hilb\times_{\tt/W}\tt]_\red$ since the scheme $\ha$ is reduced.

Further, 
Corollary \ref{P} and  Lemma \ref{irr} imply
that the scheme $\ha$ is  Cohen-Macaulay 
and smooth in codimension 1.
We conclude that $\ha$  is a normal scheme which is birational
and finite over $[\hilb\times_{\tt/W}\tt]_\red$.
This yields the isomorphism of the proposition.
The Gorenstein property of $\ha$ follows,
since $\xx$ is  Gorenstein, from
the isomorphism $\ha=(\xc\times_{\zz^\circ}\ss)/G$,
as  has been
already mentioned earlier.
\end{proof}

We remark that Haiman has proved in  \cite{Ha1} that the 
isospectral Hilbert scheme is, in fact, normal.
Assuming this result, Proposition \ref{haim}
implies that  the 
isospectral Hilbert scheme is  Cohen-Macaulay and  Gorenstein,
the main result of \cite{Ha1}. Unfortunately,
our method does not seem to yield an independent
proof of  normality of  the 
isospectral Hilbert scheme, while the proof of  normality given
in \cite{Ha1} is based on the `polygraph theorem', a key
technical result of \cite{Ha1}.

As a consequence of Haiman's work, we deduce
\begin{cor}\label{procesi} The vector bundle
$\pp$, on $\hilb$, is isomorphic to the
{\em Procesi bundle}, cf. \cite{Ha1}.
\qed
\end{cor}

\section{Geometric construction of the \hc module}\label{int}

\subsection{Filtered $\dd$-modules}\label{filt_sec}
Following G. Laumon and M. Saito, we consider
the category $\Lmod{F\dd_X}$, of {\em left} $\dd_X$-modules $M$ 
 equipped with a good filtration $F$.
This is an exact (not abelian) category
and there is an associated  derived category.
Let $D^b_{coh}(\Lmod{F\dd_X})$
be the full triangulated subcategory of 
the  derived category
whose objects are isomorphic to
bounded complexes $(M,F),$ of filtered $\dd_X$-modules, 
such that each cohomology group $\H^i(M,F)\in\Lmod{F\dd_X}$
 is a coherent
$\dd_X$-module.
The assignment
$(M,F)\mto \ggr^F M$ gives a functor
$\Lmod{F\dd_X}\to \Lmod{\oo_{T^*X}}$
that can be extended to a triangulated functor
$\dgr:\ 
D^b_{coh}(\Lmod{F\dd_X})$ $\to \dcoh( X),$
cf. \cite{Sa}, \cite{La}.

Associated with a {\em proper} morphism $f:\ X\to Y$, of
smooth varieties, there is a natural
direct image functor $\int_f$, resp. $\int_f^R$, 
between  bounded derived categories of
coherent left, resp. right,
$\dd$-modules on $X$ and $Y$, respectively, cf. \cite{Bo}. 
For  a left $\dd_X$-module $M$,
one has
$\int_f M= K\inv_Y\bo_{\oo_Y}\int_f^R (K_X\bo_{\oo_X}M).$

The direct image functor can be upgraded to
a functor
$D^b_{coh}(\Lmod{F\dd_X})\to
D^b_{coh}(\Lmod{F\dd_Y}),$
$(M,F)\mto \int_f(M,F)$,
between filtered derived categories, cf. \cite{Sa}, \cite{La}. 
The latter functor commutes with the
associated graded functor $\dgr(-)$.
We will only need a special
case of this result for  maps of the form
$f: X\times Y\to Y,$ the  projection along
a proper variety $Y$. 
In this case, one has a diagram
\beq{fdiag}
\xymatrix{
T^*Y\ &\ X\times T^*Y\ \ar@{->>}[l]_<>(0.5){\pr}
\ar@{^{(}->}[rrr]^<>(0.5){\imath_{X\times T^*Y\to T^*(X\times Y)}}
&&&\ T^*X\times T^*Y=T^*(X\times Y).
}
\eeq
Here,
$\imath: X\into T^*X$ denotes the zero section and we
put $\imath_{X\times T^*Y\to T^*(X\times Y)}:=
\imath\boxtimes\Id_{T^*Y}$. 

The relation between the functors $\dgr(-)$ and $\int_f$
is provided by the following result
\begin{thm}[\mbox{\cite{La}, \cite{Sa}}]
\label{fgr} Let $f: X\times Y\to Y$ be the
second projection where
$X$ is proper.
Then, for any $(M,F)\in D^b_{coh}(\Lmod{F\dd_{X\times Y}})$,
in $\dcoh((T^*Y))$, there is a functorial
isomorphism
$$\dgr\mbox{$\left(\int_f (M,F)\right)$}=
R\pr_\idot\big((K_X\boxtimes \oo_{T^*Y})\ 
\lo_{\oo_{X\times T^*Y}}
L\imath_{X\times T^*Y\to T^*(X\times Y)}^*
\dgr(M,F)\big).
\eqno\Box$$
\end{thm}

Let $f:\ X\to Y$ be a proper morphism
and  $(M,F)$ a filtered $\dd_X$-module.
Each cohomology 
 of the filtered complex $\int_f(M,F)$
is a $\dd_Y$-module $\H^\hdot\big(\int_f(M,F)\big)$
that comes equipped 
with an induced filtration. However,
Theorem \ref{fgr} is not sufficient, in general, for
describing $\ggr\H^\hdot\big(\int_f(M,F)\big)$,
the associated graded sheaves.
To explain this let, more generally,
$$E:\ \ldots\stackrel{d_{k-2}}\too E^{k-1}\stackrel{d_{k-1}}
\too E^k\too E^{k+1}\stackrel{d_{k}}\too\ldots,$$
be an arbitrary filtered complex in an abelian category.
Then, one has an induced
filtration on each cohomology group
$H^k(E),\ k\in\Z$, and there is a spectral
sequence $H^\hdot(\gr E)\ \Rightarrow\ \gr  H^\hdot(E)$.
The filtered complex $E$ is said to be {\em strict}
if the  morphism
$d_k: F_jE^k\to \im d_{k}\cap F_j E^{k+1}$ is surjective,
for any $k,j\in\Z$. For a strict filtered complex,
the spectral sequence  degenerates
and
one has a canonical isomorphism
$\gr  H^\hdot(E)\cong H^\hdot(\gr E)$.

Now, let $f: X\to Y$ be a projective morphism of smooth
varieties. A polarized \mhm  on $X$
may be viewed as an object $ (M,F)\in D^b_{coh}(\Lmod{F\dd_X})$.
The corresponding 
 direct image,  $\int_f(M,F)\in  D^b_{coh}(\Lmod{F\dd_Y})$,
may be thought of  as a filtered bounded complex
 of $\dd_Y$-modules.
One of the main results of Saito's theory 
reads,
see \cite[Theorem 5.3.2]{Sa}:
\begin{thm}\label{saito} 
For any
 $(M,F)\in \mh(X)$ and any projective morphism
  $f: X\to Y$, the filtered complex $\int_f(M,F)$
is {\em strict} and each cohomology group
$\H^j(\int_f(M,F))$ has a natural structure of a polarized
\mhm  on ~$Y$.
\end{thm}

In the situation of the theorem, we refer to
 the induced filtration on 
 $\H^j(\int_f(M,F)),\ j=0,1,\ldots$,
as the Hodge filtration
and let $\ggr^\hodge\H^j(\int_f(M,F))$
denote the associated graded coherent sheaf on $T^*Y$.
Similar notation will be used for right $\dd$-modules.

\subsection{The Hotta-Kashiwara construction}\label{hk}
The sheaf $\oo_\tg$ has an obvious  structure of
a holonomic left $\dd_\tg$-module. So, 
one has $\int_{\mu\times\nu}\oo_\tg$, 
 the direct image of this $\dd_\tg$-module 
via the proper morphism $\mu\times\nu$.
Each cohomology group
$\H^k(\int_{\mu\times\nu}\oo_\tg)$ is a holonomic
$\dd_{\g\times\t}$-module  set theoretically
supported on the variety $\yy\sset \g\times\t$.

 Hotta and Kashiwara proved the following important result,
\cite{HK1}, Theorem 4.2.

\begin{thm}\label{pimm} For  any
$k\neq 0$, we have
$\H^k(\int_{\mu\times\nu}\oo_\tg)=0,$
and there is a natural
isomorphism of $\DD$-modules
$\H^0(\int_{\mu\times\nu}\oo_\tg)
\cong
\DD/\I'$,
where $\I'\sset\DD$ is the following left ideal
\beq{I}
\I':=\DD\cdot(\ad\g\o 1)+\DD\cdot\{P-\rad(P)\mid P\in \C[\g]^G\}+
\DD\cdot\{Q-\rad(Q)\mid Q\in(\sym\g)^G\}.
\eeq
\end{thm}

An alternative, more direct, proof of
a variant of Theorem \ref{pimm} may be found in \cite{HK2}.

It is clear, by comparing
formulas \eqref{M} and \eqref{I}, that, using the
notation of Lemma \ref{supp},
one has an inclusion $\I'\sset\I$, of left ideals.

\begin{lem}\label{I=I} We have $\I'=\I$.
\end{lem}
\begin{proof} According to \cite[Lemma 4.6.1]{HK1},
one has that $\DD/\I'$ is either a simple $\DD$-module
or else $\DD/\I'=0$. 
We conclude, since $\DD/\I=\mm\neq 0$ by Lemma \ref{supp}(ii),
 that the natural surjection
$\DD/\I'\onto \DD/\I=\mm$ must be an isomorphism.
\end{proof}

The canonical section
$\omega$, cf. Proposition
\ref{tb}(iii), gives an element in
$\H^0(\int_{\mu\times\nu}^R K_\tg)$.
Therefore, writing  $dx$, resp. $dt$, for
 constant volume forms on
$\g$, resp. $\t$,
one can view  $(dx\o dt)\inv\o\om$ as a section of 
$\H^0(\int_{\mu\times\nu}\oo_\tg)$.
The proof of the isomorphism 
$\H^0(\int_{\mu\times\nu}\oo_\tg)\cong\DD/\I'$,
see \cite{HK1}, \S 4.7,
combined with  Lemma \ref{I=I}, implies the following

\begin{cor}\label{sec}
The assignment $u\mto u[(dx\o
dt)\inv\o\om]$ yields an isomorphism
$\mm\iso \H^0(\int_{\mu\times\nu}\oo_\tg).$
\end{cor}

\begin{rem}
It follows from the corollary
that the section $(dx\o dt)\inv\o\omega$ is annihilated
by all differential operators of the form
$u\o1-1\o\rad u,\ u\in\dd(\g)^G$.
This is a strengthening of \cite{HK1},
formulas (4.7.2)-(4.7.3).
\end{rem}

\subsection{}\label{D} 
It will be convenient to factor the map $\mu\boxtimes\nu:\ \tg\to \g\times\t$
as a composition of a closed imbedding
$\eps:\ \tg\into \bb\times\g\times\t,\ (\b,x)\mto (\b, x, \nu(x))$
and a  projection $f:\ \bb\times\g\times\t\to \g\times\t,$
along the first factor. 
The image of the imbedding $\eps$ is a smooth subvariety in
$\eps(\tg)\sset \bb\times\g\times\t$. 
Therefore, $E:=\int_\eps^R K_\tg$ is a simple
holonomic right $\dd$-module on $\bb\times\g\times\t$.
This $\dd$-module is
generated by the section $\eps_*\om$.
So, we have a surjective morphism
$\gamma:\ \dd_{\bb\times\g\times\t}\onto E,\, 1\mto \eps_*\om$,
of right  $\dd_{\bb\times\g\times\t}$-modules.
Applying the functor $\int_f^R$ one obtains a
morphism $\int_f^R\gamma:\ \int_f^R\dd_{\bb\times\g\times\t}\onto \int_f^RE$,
of complexes of right $\dd_{\g\times\t}$-modules.

\begin{lem}\label{intd} All nonzero cohomology
sheaves of the complex $\int_f^R\dd_{\bb\times\g\times\t}$
vanish and one has
an isomorphism $\H^0(\int_f^R\dd_{\bb\times\g\times\t})\cong
\dd_{\g\times\t}$, of right $\dd_{\g\times\t}$-modules.
\end{lem}
\proof
Let $\T$ denote the tangent sheaf on $\bb$.
There is a  standard Koszul
complex ${\mathcal K}^\hdot$,
with terms ${\mathcal K}^j=\dd_\bb\o_{\oo_\bb}\wedge^j \T$,
that gives a resolution of the structure
sheaf $\oo_\bb$, cf. \S\ref{dima_sec}. 
Thus, using that $H^0(\bb,\oo_\bb)=\C$ and
$H^k(\bb,\oo_\bb)=0$ for any $k\neq0$, 
 by the definition of direct image, we get
$$\int_f^R\dd_{\bb\times\g\times\t}=
R\Ga(\bb, \ {\mathcal K}^\hdot)\,\o_\C\,
\dd_{\g\times\t}=
R\Ga(\bb, \oo_\bb)\,\o_\C\,
\dd_{\g\times\t}=\C\,\o_\C\,
\dd_{\g\times\t}=\dd_{\g\times\t}.\eqno\Box
$$
\medskip

We have $\int_{\mu\times\nu}\oo_\tg=
\int_f\br{\int_\eps\oo_\tg}=
K\inv_{\g\times\t}
\bo_{\oo_{\g\times\t}}\int^R_f E.$
Hence, from Theorem \ref{pimm} and Corollary \ref{sec},
we deduce that
$\H^j\br{\int^R_f E}=0$ for all $j\neq 0$
and, moreover, one has
an isomorphism
$K_{\g\times\t}\inv\o\H^0\br{\int^R_f E}\cong\mm$.

We identify the sheaf $\dd_{\g\times\t}$
with $\dd_{\g\times\t}^{op}$ via the
trivialization of the canonical
bundle on $\g\times\t$ provided by the section  $dx\o dt$.
Let $\bar\gamma:=\H^0(\int^R_f\gamma)$ denote the map of $0$-th cohomology
groups induced by the morphism $\int^R_f\gamma$.
Thus, using Lemma \ref{intd},  we obtain
a chain of morphisms of {\em left} $\dd_{\g\times\t}$-modules
\beq{j}
\xymatrix{
\dd_{\g\times\t}=
K\inv_{\g\times\t}\o\H^0(\int_f^R\dd_{\bb\times\g\times\t})
\ar[r]^<>(0.5){\bar\gamma}&
K\inv_{\g\times\t}\o\H^0(\int^R_f E)
=\H^0(\int_{\mu\times\nu}\oo_\tg)=\mm.
}
\eeq
It is straightforward to see,
using the explicit isomorphism of Corollary \ref{sec}, that the composite
morphism in \eqref{j}
is nothing but the natural
projection $\dd_{\g\times\t}\onto \dd_{\g\times\t}/\I=\mm$.

\subsection{}\label{filtmmG}
Let  $\La$
be the total space of the
conormal bundle on $\eps(\tg)\sset \bb\times\g\times\t$.
 Thus,
 $\La$
is a smooth Lagrangian cone
in $T^*(\bb\times\g\times\t)=T^*\bb\times\gg\times\tt$.
We will view the structure sheaf $\oo_\La$
as a coherent sheaf on $T^*\bb\times\gg\times\tt$
supported on $\La$.

We consider diagram \eqref{fdiag} in the case where
$X=\bb$ and $Y=\g\times\t$.

\begin{thm}\label{cohiso} All nonzero cohomology
sheaves of the object
$R\pr_*L\imath_{\bb\times\gg\times\tt\to T^*\bb\times\gg\times\tt}^*
\,\oo_\La\in \dcoh(\gg\times\tt)$ vanish and, we have
$$\H^0\br{R\pr_*L\imath_{\bb\times\gg\times\tt\to T^*\bb\times\gg\times\tt}^*
\,\oo_\La}=\ggr^\hodge\mm.$$
\end{thm}
\begin{proof}
In section \ref{D},
we have introduced the right $\dd_{\g\times\t}$-module $E=\int^R_\eps K_\tg$.
This module 
comes equipped with the natural structure
of a polarized Hodge module. One
has $\ggr^\hodge E=q^*_\La(\eps_\idot K_\tg),$
where $q_\La:\La\into T^*(\bb\times\g\times\t)\onto \bb\times\g\times\t$
denotes the composite map. 
Let
$M:=K_{\bb\times\g\times\t}\inv\o E$,
be the corresponding left $\dd_{\bb\times\g\times\t}$-module.
Since the canonical bundle on $\tg$, resp.
on $\g\times\t$, is  trivial we find
$$K_\bb\o
Li^*(\ggr^\hodge M)=K_\bb\o
Li^* q^*_\La(\eps_\idot K_\tg\o
K_{\bb\times\g\times\t}\inv)=K_\bb\o
Li^* q^*_\La(K_\bb\boxtimes\oo_{\g\times\t})=
Li^*\oo_\La,
$$
where we have used simplified notation
$K_\bb\o(-)=(K_\bb\boxtimes \oo_{\gg\times\tt})\ 
\bo_{\oo_{\bb\times \gg\times\tt}}(-)$
and $i:=\imath_{\bb\times \gg\times\tt\to T^*\bb\times \gg\times\tt}.$
From the above isomorphisms, applying Theorem \ref{fgr}
to the left $\dd_{\bb\times\g\times\t}$-module $M$,
we get $\dgr(\int_fM)=
R\pr_\idot Li^*\oo_\La.$
Thus, since $\int_{\mu\times\nu}\oo_\tg=
\int_f(\int_\eps\oo_\tg)=\int_fM$,
we deduce
$\dgr\big(\int_{\mu\times\nu}\oo_\tg\big)=
R\pr_\idot Li^*\oo_\La.$

Now, the sheaf $\oo_\tg$ has the natural structure
of a polarized Hodge left $\dd_\tg$-module.
Thanks to Theorem  \ref{saito},
we obtain an isomorphism
$\ggr^\hodge\H^j(\int_{\mu\times\nu}\oo_\tg)
=\H^j(R\pr_\idot Li^*\oo_\La)$,
for any $j\in\Z.$
Theorem \ref{pimm} completes the proof.
\end{proof}

\section{Completing the proof of Theorem \ref{mthm}}
\subsection{Symplectic geometry interpretation}\label{sympl}
Let $\Phi:\ \wt\N\iso T^*\bb$ be an
isomorphism obtained by composing the natural
$G$-equivariant vector bundle isomorphism $\wt\N\cong T^*\bb$,
used in \S\ref{dima_sec}, with the sign involution
along the fibers of the vector bundle ~$T^*\bb$.

Recall the setup and notation of \S\ref{tx}.
The following result gives
an interpretation of the variety $\tgg$
in terms of symplectic geometry. 

\begin{prop}\label{conorm}
The  map  $\Psi=
(\Phi\ccirc\kkap)\times\mmu\times\nnu:\ \tgg\to
T^*\bb\times\gg\times\tt$ gives
an isomorphism of $\tgg$ with the variety
$\La\sset T^*\bb \times \gg\times\tt,\ $
the total space of the conormal
bundle to the subvariety $\eps(\tg).$
\end{prop}

\begin{proof} 
Fix $\b\in\bb$ and $x\in\b$.
So,
$(\b,x)\in\tg$
and the corresponding point in $\bb\times\g\times\h$
is given by the triple $ u=(\b,\,x,\,x\mo)$.
The fiber  of the vector
bundle $\TT(\bb\times\g\times\h)$
 at the point
$ u$ may be identified vith the vector
space $(\g/\b)\times\g\times\h$.
It is easy to verify that
the fiber of the sub vector bundle $\TT(\tg)$  at $ u$
is equal to the following subspace
$$\big\{(\alpha\,\operatorname{mod}\b,\,[\alpha,x]+\beta,\,\beta\mo)\
\in\ 
(\g/\b)\times\g\times\h\enspace\big|\en \alpha\in\g,\ \beta\in\b\big\}.
$$

Now, write $\langle-,-\rangle$ for an invariant bilinear form on $\g$
and use it to identify the fiber  of the vector
bundle  $\TT^*(\bb\times\g\times\h)$ at 
$ u$  vith the vector
space $[\b,\b]\times\g\times\h$.
Let $( n, y,h)\in [\b,\b]\times\g\times\h$
be a point of that vector space.
Such a point belongs to
 the fiber of $\La$ over $u$
if and only if the following equation holds
\beq{eq}
\langle\alpha, n\rangle+\langle[\alpha,x]+\beta,\, y\rangle+\langle\beta\mo,\ h\rangle=0
\quad\forall \alpha\in\g,\,\beta\in\b.
\eeq

Taking $\alpha=0$ and applying 
equation \eqref{eq}  we get $\langle[\b,\b],\, y\rangle=0$.
Hence,
$ y\in\b$ and  $h=\beta\mo$.
Next, for any $\alpha\in\g$,
we have $\langle[\alpha,x],\, y\rangle=\langle\alpha,\,[x, y]\rangle$.
Hence, for $\beta=0$ and any $ n\in[\b,\b],\  y\in\b$,
equation \eqref{eq} gives
$$0=\langle\alpha, n\rangle+\langle[\alpha,x], y\rangle=\langle
\alpha,\, n+[x, y]\rangle\quad\forall\alpha\in\g.$$
It follows that $ n+[x, y]=0$. We conclude
 that the fiber of $\La$
 over $ u$ is
equal to 
\beq{qq}\{(-[x,y],\,y,\,y\mo)\in[\b,\b]\times\g\times\h,\quad y\in\b\}.
\eeq

We have a projection
$\pr_{12}:\ \tgg\to\tg,\ (\b,x,y)\to (\b,x)$, along the last factor.
We see that the vector space \eqref{qq} is nothing but
the image of set $\pr_{12}\inv(\b,x)$ under the map
$\Psi$.
This proves that the latter map gives
an isomorphism $\tgg\iso\La$
as well as the commutativity of the diagram of the proposition.
\end{proof}

Let $pr_{\La\to T^*\bb}:\ \La\to
T^*\bb$, resp. $pr_{\La\to\gg\times\tt} 
:\ \La\to\gg\times\tt,$
denote the restriction to $\La$ of
the projection of $T^*\bb\times\gg\times\tt$
to the first, resp. along the first, factor.
By the definition $\Psi=
(\Phi\ccirc\kkap)\times\mmu\times\nnu$,
of the map $\Psi$, one has a commutative diagram
\beq{psi}
\xymatrix{
\gg\times\tt\ \ar@{=}[d]^<>(0.5){\Id}&&
\ \tgg\ \ar@{=}[d]^<>(0.5){\Psi}
\ar[rr]^<>(0.5){\kkap}\ar[ll]_<>(0.5){\mmu\times\nnu}
&&\ \wt\N\ \ar@{=}[d]^<>(0.5){\Phi}\\
\gg\times\tt\ &&\ \La\ \ar[rr]^<>(0.5){pr_{\La\to T^*\bb}}
\ar[ll]_<>(0.5){pr_{\La\to\gg\times\tt}}&&
\ T^*\bb.
}
\eeq

From this diagram,
we deduce, by Proposition \ref{conorm},
 that the map $\Psi$ induces
an isomorphism 
between the scheme $\tx=\kkap\inv\big(\imath(\bb)\big)$,
see \S\ref{txform},
and $\La\,\cap\,(\bb\times \gg\times\tt)=
pr_{\La\to T^*\bb}\inv\big(\imath(\bb)\big).$

\subsection{}\label{ulem}
Recall the open
 set $\zz^{rr}\sset \zz$, see Definition \ref{U}.
Thus, $\wt\xx^{rr}:=\mmu\inv(\zz^{rr})$ is a Zariski 
open subset of $\tx$, resp.
$\xx^{rr}$ is a  Zariski open subset of
$\zz\times_{\zz/\!/G}\hh$, and we have

\begin{lem}\label{codim1} \vi The the differential
of the morphism $\kkap: \tgg\to \wt\N$
is surjective at any point
 $(\b,x,y)\in\wt\xx^{rr};$ in particular,
the set $\wt\xx^{rr}$ is contained in the
smooth locus of the scheme $\tx$.

\vii The set $\xx^{rr}$ is contained in the
smooth locus of the variety $\xx$.
\end{lem}
\begin{proof}
Let $x\in\b$ and write $\ad_\b x$ for the map
$\b\mto [\b,\b], u\mto
[x,u]$. We  have $\dim\b-\dim\g_x\leq\dim\b-\dim\ker(\ad_\b x)=
\dim\im(\ad_\b x)\leq\dim[\b,\b]$.
If $x\in\b\cap\g^r,$ then
we have $\dim\g_x=\dim\t$ and from the above inequalities
we deduce $\dim[\b,\b]=\dim\b-\dim\g_x\leq\dim(\im\ad_\b x)
\leq\dim[\b,\b]$. It follows that one has
$\im(\ad_\b x)=[\b,\b]$.
Thus,  for $x\in\b\cap\g^r,$ the map $\ad_\b x$ is surjective.

Next, let  $y\in\b$ be another element.
The differential
of the commutator map $\kap:\
\b\times\b\to[\b,\b]$ at the point $(x,y)\in \b\times\b$
is a linear map $d_{\b,x,y}\kap:\ \b\oplus\b\to[\b,\b]$
given by the formula
$d_{\b,x,y}\kap:\ (u,v)\mto \ad_\b x(u)+\ad_\b y(v)$.
In particular, we see that
$\im(\ad_\b x)\subseteq \im(d_{\b,x,y}\kap).$

Now, let $(\b,x,y)\in\wt\xx^{rr}.$ 
Without loss of generality,
we may assume that  $x$ is a regular element of $\g$.
Thus, the map  $\ad_\b x$ is surjective.
By the preceding paragraph,
we deduce that the
map $d_{\b,x,y}\kap$ is surjective as well.
Part (i) follows from this by
$G$-equivariance.

To prove (ii),  we
use Lemma \ref{conormal}
to identify $\xx$ with $\overline{N_{\x^r}}$.
Let $q: \overline{N_{\x^r}}\to \x$ denote the
projection.
It is clear that the set $q\inv(\x^r)=N_{\x^r}$ is
contained in the smooth locus
of $\xx$.

To complete the proof
of the lemma, write $\zz^{rr}=\zz^{rr}_1\cup \zz^{rr}_2$
where $\zz^{rr}_i$ is the set of pairs $(x_1,x_2)\in\zz$
such that $x_i\in\g^r$.
It is immediate
from Proposition \ref{tb}(i)
that one has $p_{_\tt}\inv(\zz^{rr}_1)\sset q\inv(\x^r)$.
We conclude that the set $p_{_\tt}\inv(\zz^{rr}_1)$
is contained in the smooth locus of $\xx$.
By symmetry, the same holds for $p_{_\tt}\inv(\zz^{rr}_2)$
as well.
\end{proof}

We may use Proposition \ref{conorm}
to identify the set $\wt\xx^{rr}\sset\tx\sset\tgg$
with a subset of the
set $\La\cap(\bb\times\gg\times\tt)\sset
T^*\bb\times\gg\times\tt$.
Let $\ppi$ denote the restriction of the
proper morphism $\mmu\times\nnu$ to the
subset $\wt\xx^{rr}$.

\begin{cor}\label{transv} \vi The intersection
$\La\cap(\bb\times\gg\times\tt)$ is transverse
at any point of the set $\wt\xx^{rr}$.

\vii  The  map
$\ppi$ 
induces an isomorphism $\ppi:\ \wt\xx^{rr}\iso \xx^{rr}.$
\end{cor}
\begin{proof} Part (i) is equivalent to the statement of 
Lemma \ref{codim1}(i). 
To prove (ii), let $\UU$ denote the
preimage of $\zz^{rr}$ under
the first projection $\zz\times_{\zz/\!/G} \tt\to\zz$.
It is clear
that the  morphism
$\ppi$  maps $\wt\xx^{rr}$
to $\UU$.
We claim
that 
the resulting map 
$\ppi:\ \wt\xx^{rr}\to\UU$
is a set theoretic
bijection.
Indeed, it 
is surjective,
since the image of this map
contains the set $\zz^{rs}\times_{\zz/\!/G} \tt$
and $\ppi$ is a proper morphism.
To prove injectivity, we interpret the
assignment $(\b,x,y)\mto (x,y)$ 
as  the map $\pi\times\Id:\
\tg\times\g\to\yy\times\g$.
This last map gives a bijection
beteen the set $\yy^r\times\g$ and its preimage
in $\tg\times\g$, thanks to Proposition
\ref{tb}(i). Our claim follows.

Next, we observe that
since $\wt\xx^{rr}$ is a smooth scheme by 
 Lemma \ref{codim1}(i) 
the  scheme theoretic image of $\wt\xx^{rr}$
under the morphism $\ppi$ 
is actually contained 
in $\xx^{rr}=\UU_\red$.
The reduced scheme $\xx^{rr}$ is
smooth by Lemma \ref{U}.
Thus, we have a  proper morphism
$\ppi: \wt\xx^{rr}\to \xx^{rr}$,
between smooth varieties,
which is a set theoretic bijection.
Such a morphism is necessarily an isomorphism,
and part (ii) follows.

There is an alternative proof that the
morphism $\ppi: \wt\xx^{rr}\to \xx^{rr}$
is \'etale
based on symplectic geometry. In more detail,
put $X=\bb$ and $Y=\g\times\t$
and let $\eps: \tg\into X\times Y$ be 
the imbedding used in \S\ref{D}.
We have smooth locally closed subvarieties
$\x^r\sset \g\times\t$
and
$Z:=\eps(\pi\inv(\x^r))\sset X\times Y$,
respectively.
Using the notation of the
proof of Lemma \ref{codim1}(ii),
we can write $p_{_\zz}\inv(\zz^{rr}_1)\sset
q\inv(\x^r)=N_{\x^r}$. Further,
we may use Proposition \ref{conorm}
to identify $\mmu\inv(\zz^{rr}_1)$ with
an open subset of $N_Z$,
the total space of the conormal
bundle on $Z$.

We know that
the projection $X\times Y\to Y$
induces an isomorphism $Z\iso\x^r,$
by Proposition \ref{tb}(i).
Now, one can prove a general result saying 
that, in this case, the map $N_Z\to N_{\x^r}$,
induced by the projection
$T^*X\times T^*Y\to T^*Y$,
is \'etale at any point
where  $N_Z$
meets the subvariety $X \times T^*Y\sset T^*(X\times Y)$
transversely.

To complete the proof we observe
that
this transversality condition
holds for any point of the set 
  $\mmu\inv(\zz^{rr}_1)\sset N_Z,$
thanks to part (i) of the Corollary.
\end{proof}

\subsection{}\label{xyzsec}
Let $A, C,C'$ be a triple of commutative algebras.
Given a pair
$A\to C$ and $A\to C'$, of algebra homorphisms,
 one may view $C$ and $C'$ as $A$-algebras.
Therefore, the tensor product $C\o_A C'$ has the natural structure
of a commutative $A$-algebra.
One can also form $C\loo_A C'$,
a derived tensor product. The latter may be
viewed  as a commutative
DG algebra which is concentrated in nonpositive degrees
and has differential of degree
$+1$. This DG algebra is well defined up
to homotopy, cf. eg. \cite{Hi}
for details. By definition, one has
$H^{-\hdot}(C\loo_A C')=\Tor^A_\idot(C,C'),$
where we use the cohomological, rather than
homological, grading on the left hand
side. The Tor-group above is a graded
commutative algebra with degree zero
component being equal to $C\o_A C'.$

The above construction localizes. Thus, given
a smooth variety $X$ and a pair
$\cc_X$ and $\cc'_X$,  of coherent $\oo_X$-algebras, one has
a sheaf $\cc_X\loo_{\oo_X}\cc'_X$,
of DG $\oo_X$-algebras,
well defined up to quasi-isomorphism.

Now,  let $i_Y: Y\into X$,
resp. $i_Z: Z\into X$,
be closed imbeddings of smooth
subvarieties.
It is known that, in $\dcoh(X)$, there are natural isomorphisms
\beq{xyz}(i_Y)_*Li_Y^*[(i_Z)_*\oo_Z]\ \cong\
(i_Y)_*\oo_Y\lo_{\oo_X}(i_Z)_*\oo_Z\
\cong\ (i_Z)_*Li_Z^*[(i_Y)_*\oo_Y].
\eeq

We  apply the discussion above
in the case where $\cc_X=(i_Y)_*\oo_Y$
and $\cc_X'=(i_Z)_*\oo_Z$.
Thus, the object at the middle
of formula \eqref{xyz} may be viewed as
a DG $\oo_X$-algebra. There are, in fact,  DG $\oo_X$-algebra
structures on two other objects in
\eqref{xyz} as well (the DG algebra $\aa^\hdot$ 
from \S\ref{dima_sec} is an example of this).
 Furthermore, the isomorphisms in \eqref{xyz} are
DG $\oo_X$-algebra quasi isomorphisms.

We are interested
 in the special case where $X=T^*\bb\times\gg\times\tt$.
Write $pr_{T^*\bb}$,
resp. $pr_{\gg\times\tt}$
for the projection of the variety
$T^*\bb\times\times\gg\times\tt$
to the first, resp.
along the first, factor.
We will identify $\bb$ with a subvariety of $T^*\bb$
via the zero section $\imath: \bb\into T^*\bb$.

Now, in the setting of
diagram \eqref{xyz}, we
put
$Y:=\La$ and $Z:=\bb\times\gg\times\tt$.
To simplify notation, we will identify 
the structure sheaf
$\oo_Y$, resp. $\oo_Z$, with
the corresponding $\oo_X$-module 
$(i_Y)_*\oo_Y$, resp. $(i_Z)_*\oo_Z$.
Using the notation of diagram
\eqref{psi}, we can write
$pr_{\La\to T^*\bb}=pr_{T^*\bb}\ccirc
i_\La$.
We deduce a chain of DG algebra quasi-isomorphisms
$$
Li^*_\La(\oo_{\bb\times\gg\times\tt})=
Li^*_\La[pr_{T^*\bb}^*(\imath_*\oo_\bb)]
=L(pr_{T^*\bb}\ccirc i_\La)^*(\imath_*\oo_\bb)
=Lpr_{\La\to T^*\bb}^*(\imath_*\oo_\bb).
$$

Next, we use the isomorphism
$\Phi: \wt\N\iso T^*\bb$
to identify
the zero section $\imath: \bb\into \wt\N$
with the zero section $\imath: \bb\into T^*\bb$.
Thus, from the above chain
of isomorphisms, by commutativity of diagram \eqref{psi}, we get
\beq{xy}
R(pr_{\La\to \gg\times\tt})_* Li_Y^*(\oo_Z)=
R(pr_{\La\to \gg\times\tt})_*
Lpr_{\La\to T^*\bb}^*(\imath_*\oo_\bb)\cong
R(\mmu\times\nnu)_*L\kkap^*(\imath_*\oo_\bb).
\eeq

Using the notation
of Theorem \ref{cohiso}, we can write
$i_Z=\imath_{\bb\times\times\gg\times\tt
\to T^*\bb\times\gg\times\tt}$, resp. $pr_{\gg\times\tt}=\pr$.
Now, we apply the functor $(Rpr_{\gg\times\tt})_*(-)$ to each of
the 3
objects of our diagram \eqref{xyz}. Combining  the isomorphisms
in  \eqref{xyz} with \eqref{xy}, we
obtain the following  DG algebra quasi-isomorphisms
\beq{dg1}
R\mmu_*L\kkap^*(\imath_*\oo_\bb)
\cong
R(pr_{\gg\times\tt})_*(\oo_{\bb\times\gg\times\tt}\!
\underset{\oo_{T^*\bb\times\gg\times\tt}}\lo
\!\oo_\La)
\cong
R\pr_*L\imath_{\bb\times\times\gg\times\tt
\to T^*\bb\times\gg\times\tt}^*\oo_\La.
\eeq
 
Note that
$L\kkap^*(\imath_*\oo_\bb)=\aa^\hdot$ is the DG algebra considered
in \S\ref{dima_sec}. Thus,
Theorem \ref{dima2} follows from the above isomorphisms and Theorem
\ref{cohiso},
cf. also \S\ref{keypf}.

\subsection{Proof of Lemma \ref{ordhodge}}\label{pf}
First of all, we specify our normalization
of the Hodge filtration on the \hc module $\mm$.
To this end, we recall the setting of \S\ref{D}.
Thus, we have a simple holonomic right
$\dd_{\bb\times\g\times\t}$-module $E$, with support
$\eps(\tg)$, and a surjective
homomorphism $\gamma:\ \dd_{\bb\times\g\times\t}
\onto E$, of right $\dd_{\bb\times\g\times\t}$-modules.
The order filtration on
$\dd_{\bb\times\g\times\t}$ induces,
via the projection $\gamma$, a quotient
filtration on $E$. The resulting filtration
is known to be equal to the Hodge filtration on $E$
up to a shift, since $\eps(\tg)$ is a smooth
submanifold in $\bb\times\g\times\t$.
We choose the normalization of
the Hodge filtration on $E$ so that
the two filtrations  coincide.

The filtration on $E$, resp. on $\dd_{\bb\times\g\times\t}$,
makes $\int^R_f E$, resp. $\int^R_f\dd_{\bb\times\g\times\t}$,
a filtered complex.
We define the Hodge filtration
on $\mm$ to be the filtration
that comes from the  induced filtration on $\H^0(\int^R_f E)$
via the isomorphism $\mm=
K_{\g\times\t}\inv\o\H^0(\int^R_f E)$, cf. \eqref{j}.

\proof[Proof of part (i)] With the above choice of normalization,
 the map $\gamma$
becomes a morphism of filtered $\dd$-modules.
It follows that the induced morphism
$\int^R_f\gamma: \int^R_f\dd_{\bb\times\g\times\t}\to
\int^R_f E$ is a morphism of filtered complexes.

The proof of Lemma \ref{intd} shows
that
the  filtration on the $\dd$-module
$\H^0(\int^R_f\dd_{\bb\times\g\times\t})$,
induced by the filtered structure on $\int^R_f\dd_{\bb\times\g\times\t}$,
goes, under
the isomorphism of the Lemma,
to the standard order filtration on
the sheaf $\dd_{\g\times\t}$.
It follows that all the maps in 
\eqref{j} are filtration preserving.
Thus, writing $\bar\gamma: \dd_{\g\times\t}
\onto \mm$ for the
composite map in \eqref{j},
we get $\bar\gamma(F_k^\ord\dd_{\g\times\t})
\sset F^\hodge_k\mm$, for any $k\in\Z$.
Now, according to the remark after
formula \eqref{j}, the map $\bar\gamma$ is  the
natural
projection $\dd_{\g\times\t}\onto
\dd_{\g\times\t}/\I=\mm$.
This yields the inclusions of part (i) of Lemma \ref{ordhodge},
since the order filtration $F^\ord\mm$,
on the \hc module, was defined as the
quotient filtration on $\dd_{\g\times\t}/\I$.

\proof[Proof of part (ii)] 
Recall the setup of
the proof of Theorem \ref{cohiso} and
introduce simplified notation $i:=\imath_{\bb\times\gg\times\tt\to
T^*\bb\times\gg\times\tt}.$
The DG algebra structure on the object on the right
of  \eqref{dg1}
induces an $\oo_{\gg\times\tt}$-algebra structure
on  $\H^0(R\pr_*Li^*\oo_\La),$
the 0-th cohomology sheaf.
The  $\oo_{\gg\times\tt}$-algebra structure on 
$\ggr^\hodge\mm$ referred to in part (ii) of Lemma \ref{ordhodge}
is induced by that on $\H^0\big(R\pr_*Li^*\oo_\La)$
via the isomorphism of  Theorem \ref{cohiso}.

Observe next
that we have
$R\pr_*Li^*\oo_{T^*\bb\times\gg\times\tt}=
R\pr_*(\oo_\bb\boxtimes\oo_{\gg\times\tt})=
\oo_{\gg\times\tt}$, since $H^k(\bb,\oo_\bb)=0$
for all $k>0$, cf. proof of Lemma \ref{intd}.
One may upgrade the above to a quasi-isomorphism 
$\oo_{\gg\times\tt}\iso R\pr_*Li^*\oo_{T^*\bb\times\gg\times\tt},$
of
DG algebras.
Thus, applying  the functor $R\pr_*Li^*(-)$
to an obvious restriction
morphism  $\oo_{T^*\bb\times\gg\times\tt}
\to\oo_\La$   yields   DG algebra morphisms
\beq{mor1}
\oo_{\gg\times\tt}\ \iso\
R\pr_*Li^*\oo_{T^*\bb\times\gg\times\tt}
\ \to\
R\pr_*Li^*\oo_\La.
\eeq
Applying further the functor $\H^0(-)$ 
we obtain a chain of $\oo_{\gg\times\tt}$-algebra morphisms
\beq{mor}
u:\ \oo_{\gg\times\tt}\ \iso\
\H^0(R\pr_*Li^*\oo_{T^*\bb\times\gg\times\tt})
\ \to\ 
\H^0(R\pr_*Li^*\oo_\La)\ \iso\ \ggr^\hodge\mm,
\eeq
where the last map
is an algebra morphism thanks to  our definition of the algebra structure
on $\ggr^\hodge\mm$.

Let $u$ be the composite morphism in \eqref{mor}.
Further, let $u_{\text{lem}}$ denote the composite morphism
involved in the statement of Lemma  \ref{ordhodge}(ii)
and let $\text{\it {res}}:\
\oo_{\gg\times\tt}\onto\oo_{\zz\times_{\zz/\!/G}\tt}$
be the natural restriction map.
It is not hard to see from the proof of part (i)
of the lemma (repeat  arguments
in the proof of Theorem \ref{cohiso} in
a simpler case where the $\dd$-module
$E$ is replaced by 
 $\dd_{\bb\times\g\times\t}$)
one has $u=u_{\text{lem}}\ccirc \text{\it {res}}$. 
We know that both $u$ and $\text{\it {res}}$ are algebra
maps. It follows that $u_{\text{lem}}$
is an algebra map as well.\qed

\subsection{Proof of Lemma \ref{keylem}}\label{keypf}
We choose a Zariski open subset 
$U\sset\gg\times\tt$ such that one has $\xx\cap U=\xx^{rr}$.
We claim 
that 
$\H^k(R\pr_*Li^*\oo_\La)|_U=0$
for all $k\neq 0$ and, moreover,  the restriction of the composite morphism  in
\eqref{mor1} to the set $U$ 
yields a quasi-isomorphism
\beq{mor2}
\oo_{\xx^{rr}}=\oo_{\gg\times\tt}|_{\xx^{rr}}
\ \iso\
\H^0(R\pr_*Li^*\oo_\La)|_U.
\eeq

To prove the claim, recall the setting  of formula \eqref{dg1}
and the DG algebra $\aa^\hdot$, see \S\ref{dima_sec}.
Let  $\wt U:=\mmu\inv(U)\sset\tgg$, and put $\wt\xx^{rr}:=
\mmu\inv(\xx^{rr})=
\tx\cap\wt U.$
Further, view $\aa^0=\oo_\tx$ as a DG algebra with zero
differential.
The transversality result
from part (i) of Corollary \ref{transv} 
implies that we have
$\H^k\big(L\kkap^*(\imath_*\oo_\bb)\big)|_{\wt U}=0,
\ \forall
k>0$ and, moreover, that the morphism
$\oo_{\wt U}\to L\kkap^*(\imath_*\oo_\bb)|_{\wt U}$
induced by the  DG algebra
imbedding $\aa^0|_{\wt U}\into \aa^\hdot|_{\wt U}$
descends to a quasi-isomorphism
$\oo_{\wt\xx^{rr}}\iso L\kkap^*(\imath_*\oo_\bb)|_{\wt U}.$

Let $\ppi:=(\mmu\times\nnu)|_{\wt U}.$
We apply the functor
$R\ppi_*$ to the above quasi-isomorphism
and use   Corollary \ref{transv}(ii).
We deduce that the  composite  of 
canonical
morphisms $\dis\oo_U\to R\ppi_*\oo_{\wt U}\to
R\ppi_* \big(L\kkap^*(\imath_*\oo_\bb)|_{\wt U}\big)$
descends to  a quasi-isomorphism
$\dis \oo_{\xx^{rr}}\iso
\big(R\ppi_*\ccirc L\kkap^*(\imath_*\oo_\bb)\big)|_U.$

Next, we use the isomorphism
between the objects on the left
 and on the right of formula \eqref{dg1}.
Thus, we have
$(R\pr_*Li^*\oo_\La)|_U\cong
\big(R(\mmu\times\nnu)_*L\kkap^*(\imath_*\oo_\bb)\big)|_U.$
By the conclusion of the previous paragraph, we obtain
that  the 
canonical map $w:\ \oo_U\to(R\pr_*Li^*\oo_\La)|_U$
descends to  a quasi-isomorphism
$\oo_{\xx^{rr}}\iso
(R\pr_*Li^*\oo_\La)|_U$.
The  morphism $w$ here corresponds, via the isomorphisms
from \eqref{dg1}, to 
the restriction to $U$ of
the composite morphism 
in \eqref{mor1}.
It follows that the map $u$,  in \eqref{mor},
descends to  an isomorphism
$u|_{_U}:\ \oo_{\xx^{rr}}\iso
\H^0(R\pr_*Li^*\oo_\La)|_U$.
This proves our claim  that
\eqref{mor2} is an isomorphism.

Now, the proof of Lemma \ref{ordhodge}(ii) shows
that the map $X=\Spec(\ggr^\hodge\mm)\to \gg\times\tt$
is induced by the $\oo_{\gg\times\tt}$-algebra morphism $u$,
 in \eqref{mor}. The above map factors through
the map $f: X\to \xx$. Thus, the map
$f:\ f\inv(\xx^{rr})\to \xx^{rr}$ may be identified
with the map induced by the composite
algebra morphism $u|_{_U}:\ \oo_{\xx^{rr}}\to
\H^0(R\pr_*Li^*\oo_\La)|_U\iso\ggr^\hodge\mm|_U$.
Formula \eqref{mor2} says that the
latter morphism is an isomorphism, completing the proof of
Lemma \ref{keylem}.
\qed

\section{Some technical results}\label{END}
\subsection{}\label{ss1}
We write $x=h_x+n_x$ for the Jordan decomposition of an 
element $x\in\g$.
We say that $(x,y)\in\gg$ is a semisimple pair
if both $x$ and $y$ are semisimple elements of $\g$.

\begin{lem}\label{lss} \vi 
For  $(x,y)\in\gg$, the following conditions are equivalent
\begin{itemize}

\item
 There exists a semisimple pair $(h_1, h_2)\in\tt$
such that, for any polynomial $f\in\C[\gg]^G$, we have
$f(x,y)=f(h_1,h_2)$.

\item
There exists a Borel subalgebra $\b\sset\g$ such that
$x,y\in\b$.
\end{itemize}

\vii If $(x,y)\in \zz$ then
the  $G$-diagonal orbit of the pair $(h_x,h_y)$ is the unique closed
$G$-orbit contained in the closure
of  the $G$-diagonal orbit
of $(x,y)$.
\end{lem}
\begin{proof} Let
$\t\sset \b$ be Cartan and Borel subalgebras of $\g$
and let $T$ be the
maximal torus corresponding to $\t$. 
Clearly, there exists a suitable one parameter
subgroup $\gamma:\ \C^\times\to T$
such that, for any   $t\in\t$ and $n\in [\b,\b]$,
one has
$\underset{^{z\to0}}\lim\ \Ad \gamma(z)(t+n)=t$.

To prove (i), let $(x,y)\in \gg$. We can write
$x=h_1+n_1,\ y=h_2+n_2,$ where $h_i\in\t$ and $n_i\in[\b,\b]$.
We see from the above that the pair $(h_1,h_2)$
is contained in the closure of the
$G$-orbit of the pair $(x,y)$ hence,
for any $f\in\C[\gg]^G$, we have
$f(x,y)=f(h_1,h_2)$.

Conversely, let $(x,y)\in \gg$ be such that
$f(x,y)=f(h_1,h_2)$ holds for  any $f\in\C[\gg]^G$.
The group $G_{h_1,h_2}$  is clearly a reductive subgroup
of $G$.
Hence, the $G$-diagonal orbit of $(h_1,h_2)$
is closed in $\gg$, cf. eg. \cite{GOV}.
It follows  that this $G$-orbit is the unique
closed orbit that is
contained in the closure
of the $G$-diagonal of the pair $(x,y)$,
since $G$-invariant polynomials on $\gg$
separate closed $G$-orbits.
By the Hilbert-Mumford criterion, we deduce
that  there exists a suitable one parameter
subgroup $\gamma:\ \C^\times\to G$
such that, conjugating the pair $(h_1,h_2)$
if necessary, on gets
$\underset{^{z\to0}}\lim\ \Ad \gamma(z)(x,y)=(h_1,h_2)$.

Now, let $\fa$ be the Lie subalgebra of $\g$
generated by the elements $x$ and $y$.
We deduce from the above that, for any
$u\in [\fa,\fa]$, one has
$\underset{^{z\to0}}\lim\ \Ad \gamma(z)(u)=0$.
This implies that $u$ is a nilpotent element
of $\g$. Thus, $[\fa,\fa]$ is a nilpotent Lie algebra,
by Engel's theorem. We conclude that $\fa$ is
a solvable  Lie algebra. Hence, there exists
a Borel subalgebra $\b$ such that $\fa\sset\b$,
and (i) is proved.

To prove (ii), we observe that
the elements $h_x,h_y,n_x,n_y$ generate
an abelian Lie subalgebra of $\g$.
Hence, there exists a Borel subalgebra
$\b$ such that we have $h_x,h_y,n_x,n_y\in\b$.
We may  choose a Cartan subalgebra $\t\sset\b$ such that
$h_x,h_y\in\t$.
Then, the argument at the beginning of the proof shows
that the pair $(h_x,h_y)$
is contained in the closure of the $G$-diagonal
orbit of $(x,y).$ 
\end{proof}

\subsection{}\label{ss2} Below, we will have to consider several
 reductive
Lie algebras at the same time. To avoid confusion,
we write $\zz(\fl)$ for the commuting
scheme of a  reductive
Lie algebra $\fl$, and use similar notation
for other objects associated with $\fl$.

Let  ${\mathfrak N}(\fl)$ be the 
{\em nilpotent} commuting variety of $\fl$,
the variety 
of pairs of commuting {\em nilpotent} elements
of $\fl$, equiped with reduced scheme structure.
It is clear that we have ${\mathfrak N}(\fl)=
{\mathfrak N}(\fl')\sset
\fl'\times\fl'$,
where $\fl':=[\fl,\fl]$ is the derived Lie algebra
of the reductive Lie algebra $\fl$.
According to \cite{Pr}, the irreducible components
of nilpotent commuting variety are parametrized
by the conjugacy classes of {\em distinguished}
nilpotent elements  of $\fl'$. The irredicible
component corresponding to such a conjugacy class
is equal to the closure in $\fl'\times\fl'=T^*(\fl')$
of the total space of the conormal bundle on that
conjugacy class. It follows, in particular, that
the dimension of each irredicible component of the variety
${\mathfrak N}(\fl)$ equals $\dim\fl'$.

Now, we fix a reductive connected group $G$ with Lie algebra $\g$,
and a
 Cartan subalgebra
 $\t\sset\g$.
Recall that the centralizer of an element of $\t$ is called
a standard Levi subalgebra of $\g$. Let $S$ be the set
of standard Levi subalgebras.
Given a  standard Levi subalgebra $\fl$,
let  $\t_\fl$ denote the center of $\fl$,
and let 
$\tts_\fl\sset\tt$,
 denote the set of elements
$(t_1,t_2)\in \tt$ such that $\g_{t_1,t_2}=\fl$.
This is  clearly a Zariski open
subset of $\t_\fl\times\t_\fl$.



Recall the isospectral commuting variety $\xx(\g)$
and the projection $p_{_\tt}:
\xx(\g)\to\tt.$ 
For each standard Levi subalgebra $\fl$, of
$\g$, we put $\xx_\fl(\g):=(p_{_\tt})\inv(\tts_\fl),$
and view this set
as a reduced scheme. Let  $L$ denote the Levi subgroup of $G$
such that $\fl=\Lie L$.

\begin{lem}\label{codim} For any standard Levi subalgebra $\fl$, 
the following assignment
$$G\times\gg\times\tts_\fl\ni
g\times (y_1,y_2)\times (t_1,t_2)\mto
 \big(\Ad g(y_1+t_1),\,\Ad g(y_2+t_2)\big)\times (t_1,t_2)
\in\gg\times\tts_\fl$$
induces a $G$-equivariant isomorphism
\beq{fiber}
\big(G\times_{ L}{\mathfrak N}(\fl)\big)\,\times\,\tts_\fl
\ \iso\ \xx_\fl(\g).
\eeq

The map $p_{_\tt}:\, \xx_\fl(\g)\to\tts_\fl$
goes, under the  isomorphism,
to the  second projection $\big(G\times_{ L}{\mathfrak N}(\fl)\big)\,\times\,
\tts_\fl\to\tts_\fl.$ Thus,
all irreducible
components  of the set $\xx_\fl(\g)$
have the same dimension and are in one-to-one correspondence
with the distinguished nilpotent cojugacy classes in
the Lie algebra $\fl$.
\end{lem}

\begin{proof} It is clear that any
commuting semisimple pair 
is  $G$-conjugate to a pair of elements
of $\t$.
Lemma \ref{lss}(ii) implies that,
for any $(x_1,x_2,t_1,t_2)\in\zz(\g)\times_{\zz(\g)/\!/G}\tt$
and any polynomial
$f\in\C[\zz(\g)]^G$, we have $f(x_1,x_2)=f(h_{x_1},h_{x_2})$.
We know that
  $W$-invariant polynomials separate
$W$-orbits in $\tt$ and that the
restriction map \eqref{Jo} is surjective.
We deduce
that the pair $(h_{x_1},h_{x_2})$ is
$W$-conjugate to the pair $(t_1,t_2)$.
It follows that the pair $(x_1,x_2)$ is
$G$-conjugate to a pair 
of the form $(t_1+y_1,t_2+y_2)$,
for some $(y_1,y_2)\in {\mathfrak N}_{\fl'}$.
The isomorphism of the lemma easily follows from this.
\end{proof}

We are now ready to comple the proof of Lemma
\ref{xxbasic}.

\begin{cor}\label{xdense} The set $\xx^{rs}$ is Zariski
dense in $\xx$.
\end{cor}
\begin{proof} The set $\tt^r$ is clearly dense in
$\tt$, hence, the closure of the set $\zz^{rs}(\g)$
contains $\zz^s(\g)$, the set of all semisimple pairs
in $\zz(\g)$. 
It follows that  the closure of the set $\xx^{rs}(\g)$
in $\xx$
contains the set $p_{_\zz}\inv(\zz^s(\g))$.
Thus, to prove the Corollary, it suffices to show that,
for any $w\in W$ and 
any $(t_1,t_2)\in\tts_\fl,$ the
quadruples  $(x_1,x_2,t_1,t_2)\in (p_{_\tt})\inv(t_1,t_2)$
such $(x_1,x_2)\in \zz^s(\g)$ form a dense
subset of the fiber $(p_{_\tt})\inv(t_1,t_2)$.

To prove this last statement, we use the isomorphism of
Lemma \ref{codim}. We see from the isomorphism
that it suffices to show that
${\mathfrak N}(\fl)$, the nilpotent commuting
variety of any standard Levi algebra $\fl$,
is contained in the closure of the set $\zz^s(\fl)$.
But we have $\zz^s(\fl)\supseteq\zz^{rs}(\fl)$
and 
the set $\zz^{rs}(\fl)$ is dense in $\zz(\fl)$,
by  Proposition
\ref{zzbasic}(i); explicitly, this
is Corollary 4.7 from \cite{Ri1}. The result follows.
\end{proof}

\subsection{Proof of Lemma \ref{irr}}\label{ss3} 
Let $\zz^{rr}_1$, resp. $\zz^{rr}_2$, be the set of pairs $(x_1,x_2)\in\zz_\g$
such that $x_1$, resp. $x_2$,  is a regular element of $\g$.
By Definition \ref{U}, we have $\zz^{rr}=\zz^{rr}_1\cup \zz^{rr}_2$.

\begin{lem}\label{u1}
The set $\zz^{rr}$ is a Zariski open subset
contained in the smooth locus of the scheme
$\zz$.
\end{lem}

\begin{proof}
For any $(x,y)\in \zz^{rr}_1$, we have $\dim(\g_x\cap\g_y)\leq\dim\g_x=\dim\t$,
since $x$ is a regular element.
Hence $(x,y)$ is a smooth point of $\zz$,
by Proposition \ref{zzbasic}(ii).
Further, the set $\zz^{rr}_1$ is the preimage 
of the set of regular elements under the
first projection $\zz\into\g\times\g\to\g$.
It follows that $\zz^{rr}_1$ is  open in $\zz$.
Similar arguments apply to the set $\zz^{rr}_2$,
hence, also to $\zz^{rr}$.
\end{proof}

Next, let $\fl\subsetneq\g$
be a proper standard Levi subalgebra of $\g$.
Then, $\t_\fl$, the center of  $\fl$,
has codimension 1 in $\t$ iff
$\fl$ is a {\em minimal}  Levi subalgebra of $\g$.
In that case, there is a root $\alpha\in \t^*$, in
 the root system $(\g,\t)$, such
that one has 
$\t_\fl=\ker\alpha$. This is a
 hyperplane in $\t$ and
we have $\fl=\t_\fl\oplus\sl_2$,
where
 $\sl_2=\fl'$ is the standard $\sl_2$-subalgebra
of $\g$ associated with the root $\alpha$.
Recall further that any nonzero element
$x\in\sl_2$ is regular and that
elements $x,y\in\sl_2$ commute iff
they are proportional to each other
(where the zero element is declared to be
proportional to any
element).

\begin{lem}\label{u2} We have 
 $\dim(\zz\sminus\zz^{rr})\leq\dim\zz-2.$
\end{lem}
\begin{proof}
We have the
projection
$p_{_\zz}:\xx(\g)\to\zz(\g)$ and, for
any standard Levi subalgebra $\fl\sset\g$,
put $\zz_\fl(\g)=p_{_\zz}(\xx_\fl(\g))$,
where we use the notation of Lemma \ref{codim}.
Thus, one has $\zz(\g)=\cup_{\fl\in S}\ \zz_\fl(\g).$
Note that,
for a pair of standard Levi subalgebra $\fl_1,\fl_2\sset\g$,
 the corresponding pieces $\zz_{\fl_1}$ and $\zz_{\fl_2}(\g)$
are equal whenever the Levi subalgebras
$\fl_1$ and $\fl_2$ are conjugate in $\g$;
otherwise these two pieces are disjoint.

Let $\fl$ be a   standard Levi subalgebra in $\g$.
We are interested in the dimension of the set
$\zz_\fl(\g)\sminus \zz^{rr}$.
From the isomorphism of Lemma \ref{irr}
we deduce that 
$$\dim(\zz_\fl(\g)\sminus \zz^{rr})=
\dim(G/L) + \dim\left(\big(({\mathfrak N}(\fl)\times\tts_\fl)\big)
\sminus \zz^{rr}\right).
$$
We put
$d_\fl:=\dim({\mathfrak N}(\fl)\times\tts_\fl)-
\dim\big([{\mathfrak N}(\fl)\times\tts_\fl)]
\sminus \zz^{rr}\big),$
a nonnegative integer.
Since
$\dim\zz(\g)=\dim (G/L)+\dim\fl+\dim\t$
and $\dim({\mathfrak N}(\fl)\times\tts_\fl)=
\dim\fl'+2\dim\t_\fl$,
we find
\begin{align}
\dim\zz-\dim(\zz_\fl(\g)\sminus \zz^{rr})=
\dim\fl+\dim\t-\dim\left(\big(({\mathfrak N}(\fl)\times\tts_\fl)\big)
\sminus \zz^{rr}\right)\label{ineq}\\
=\dim\fl+\dim\t+d_\fl-\dim\fl'-2\dim\t_\fl=
\dim(\t/\t_\fl) +d_\fl.\nonumber
\end{align}

To prove the lemma,
we will show that, for any standard Levi subalgebra $\fl\sset\g$,
the set $\zz_\fl(\g)\sminus \zz^{rr}$
has codimension $\geq 2$ in $\zz(\g)$.
We see from \eqref{ineq} that
this holds automatically whenever
$\dim(\t/\t_\fl)\geq 2$.
Thus, the only nontrivial cases are:
$\dim(\t/\t_\fl)=0$ or $\dim(\t/\t_\fl)=1$

In the first case we have $\t=\t_\fl$, so
 $\fl=\t$. In this case, $\tts_\fl=\tt^r$
and ${\mathfrak N}(\fl)=\{0\}$.
Thus, the set  $\tt^r\sminus \zz^{rr}$ consists of
the pairs $(h_1,h_2)\in\tt^r$ such that
neither $h_1$ nor $h_2$ is regular.
Therefore, each of these two elements
belongs to some root hyperplane in $\t$,
that is, belongs to a finite union of codimension
1 subspaces in $\t$.
We conclude that $d_\fl=\dim\tt^r-\dim(\tt^r\sminus \zz^{rr})
\leq \dim\tt^r-2$, hence, $d_\fl\geq 2$.

Next,  put $Z:={\mathfrak N}(\fl)+\tts_\fl$.
To complete the proof of the lemma,
we must show that, in 
the case where $\dim(\t/\t_\fl)=1$,
the set $Z\sminus \zz^{rr}$
has codimension $\geq 1$ in $Z$.
To this end, pick $h\in \ts_\fl$.
Then,  the centralizer of $h$ in $\g$
equals $\fl$. Therefore, for any
nilpotent element $n\in \fl^r$,
the element $h+n$ is a regular element of $\g$.
Clearly,  $(n,n)\in {\mathfrak N}(\fl)$ 
and we have
$(h+n,h+n)\in
Z\cap \zz^{rr}$.
Hence, $Z\cap \zz^{rr}$
is a nonempty Zariski open subset of $Z$.
Thus, the set
$Z\sminus \zz^{rr}$
is a closed proper subset of $Z$.
Since $Z$ is irreducible,
we conclude that the set
$Z\sminus \zz^{rr}$
has codimension $\geq 1$ in $Z$,
and we are done.
\end{proof}
\subsection{Acknowledgements} We are very grateful
to Dmitry Arinkin for several helpfull discussions.

\small{
\bibliographystyle{plain}

}

\end{document}